\documentclass[preprint,12pt]{elsarticle}

\usepackage{amssymb}
\usepackage{relsize}
\usepackage{enumerate}
\usepackage{bm}
\usepackage{amsmath}

\newtheorem{thm}{Theorem}
\newtheorem{lem}{Lemma}

\begin{document}

\begin{frontmatter}



\title{Multiscale computational method for heat conduction problems of composite structures with diverse periodic configurations in different subdomains}


\author[label1]{Hao Dong\corref{cor1}}\ead{donghaonwpu@mail.nwpu.edu.cn}
\cortext[cor1]{Corresponding author.}
\author[label2]{Junzhi Cui}
\author[label3]{Yufeng Nie}
\author[label3]{Zihao Yang}

\address[label1]{School of Mathematics and Statistics, Xidian University, Xi'an 710071, PR China}
\address[label2]{LSEC, ICMSEC, Academy of Mathematics and Systems Science, Chinese Academy of Sciences, Beijing 100190, PR China}
\address[label3]{Department of Applied Mathematics, School of Science, Northwestern Polytechnical University, Xi'an 710129, PR China}

\begin{abstract}
This study develops a novel multiscale computational method for heat conduction problems of composite structures with diverse periodic configurations in different subdomains. Firstly, the second-order two-scale (SOTS) solutions for these multiscale problems are successfully obtained based on asymptotic homogenization method. Then, the error analysis in the pointwise sense is given to illustrate the importance of developing SOTS solutions. Furthermore, the error estimates for the SOTS approximate solutions in the integral sense is presented. In addition, a SOTS numerical algorithm is proposed to effectively solve these problems based on finite element method. Finally, some numerical examples verify the feasibility
and effectiveness of the SOTS numerical algorithm we proposed.
\end{abstract}

\begin{keyword}
heat conduction problems \sep asymptotic homogenization method \sep diverse periodic configurations \sep error estimates \sep SOTS numerical algorithm
\end{keyword}

\end{frontmatter}


\section{Introduction}
\label{}
With the rapid development of material science and technology, composite materials are widely used in aeronautic and aerospace engineering owing to their excellent physical properties. These composites usually serve under complex and extreme heat environments. In practical engineering application, engineers and designers always adopt diverse composites in different subdomains of a whole engineering structure for creating the more complex components and systems. In order to obtain the optimal design of engineering structures, it is significant to accurately evaluate and predict the thermal responses of the composites. It is well-known that the complexity and heterogeneity (inclusions or holes) of composite materials often cause costly computational efforts \cite{R1,R2,R3,R4,R5,R6}. Fortunately, in the last thirty years, mathematicians and engineers have developed some multiscale methods to solve this difficult problem, such as the asymptotic homogenization method (AHM), heterogeneous multiscale method (HMM), variational multiscale method (VMS), multiscale finite element method (MsFEM) and multiscale eigenelement method (MEM), etc \cite{R7,R8}. As far as we know, many studies were performed on heat conduction problems of the composites. However, most of these studies focused on heat conduction problems of engineering structures manufactured by the same composites in the whole structures \cite{R1,R2,R3,R4,R5,R6}. To the best of our knowledge, there is a lack of adequate researches on heat conduction problems of composite structures with diverse periodic configurations in different subdomains \cite{R9,R10}.

The aim of this paper is to develop a multiscale computational method to effectively solve the heat conduction problems of composite structures with diverse periodic configurations in different subdomains. Based on the asymptotic homogenization method, we establish a novel SOTS analysis method and associated numerical algorithm for the above-mentioned multiscale problems.

This paper is outlined as follows. In Section 2, the detailed construction of the SOTS solutions for heat conduction problems of composite structures with diverse periodic configurations in different subdomains is given by multiscale asymptotic analysis. In Section 3, the error analysis in the pointwise sense of first-order two-scale (FOTS) solutions and SOTS solutions is obtained. By comparing the results of error analysis of FOTS solutions and SOTS solutions in the pointwise sense, we theoretically explain the importance of SOTS solutions in capturing micro-scale information. Moreover, an
explicit convergence rate for the SOTS solutions are derived under some hypotheses. In Section 4, a SOTS numerical algorithm based on FEM and FDM is presented to effectively solve these multiscale problems. In Section 5, some numerical results are shown to verify the validity of our SOTS algorithm. Finally, some conclusions are given in Section 6.

For convenience, throughout the paper we use the Einstein summation convention
on repeated indices.
\section{SOTS analysis of the governing equation}
Let us consider the following governing equation for heat conduction problems of composite structures with diverse periodic configurations in different subdomains
\begin{equation}
\left\{ \begin{aligned}
&-\frac{\partial }{\partial {x_i}}\Big( {k_{ij}^{\bm{\varepsilon }}({\mathbf{x}})\frac{\partial T^{\bm{\varepsilon }}(\mathbf{x})}{\partial {x_j}}} \Big)= h(\mathbf{x}),\;\;\text{in}\;\;\Omega=\bigcup\limits_{s = 1}^K{\Omega _s},\;\;\bm{\varepsilon}=(\varepsilon _1,\cdots\varepsilon_K),\\
&T^{\bm{\varepsilon}}(\mathbf{x}) = \widehat T(\mathbf{x})\;\;\;\text{on}\;\;\partial {\Omega_T},\\
&{k_{ij}^{\bm{\varepsilon }}({\mathbf{x}})\frac{\partial T^{\bm{\varepsilon }}(\mathbf{x})}{\partial {x_j}}}{n_i} = \bar q(\mathbf{x})\;\;\;\text{on}\;\;\partial {\Omega_q}
\end{aligned} \right.
\end{equation}
where $\Omega$ is a bounded convex domain in $\mathbb{R}^\mathcal{N}(\mathcal{N}=2,3)$ with a boundary $\partial\Omega=\partial\Omega_{T}\cup\partial\Omega_{q}$; $\Omega_s$ denotes a component, that is made of composite materials with periodic configuration $C_s$ and characteristic periodic size $\varepsilon_s$; $T^{\bm{\varepsilon}}(\mathbf{x})$ are undetermined temperature field; $k_{ij}^{\bm{\varepsilon }}({\mathbf{x}})$ is the second order thermal conductivity tensor; $h(\mathbf{x})$ is the internal heat source; $\widehat{T}({\bf{x}})$ is the prescribed temperature on the boundary $\partial\Omega_{T}$; $\bar{q}({\bf{x}})$ is the heat flux prescribed normal to the boundary $\partial\Omega_{q}$ with the normal vector $n_i$.

To begin with, let us set $\displaystyle\mathbf{y}^s=\frac{\mathbf{x}}{\varepsilon_s}$ for $\mathbf{x}\in \Omega_s$ as micro-scale coordinates of periodic unit cell $Q^s=(0,1)^\mathcal{N}$. With this notation, we have the following chain rule
\begin{equation}
\displaystyle\frac{\partial}{\partial x_i}\rightarrow\frac{\partial}{\partial x_i}+\frac{1}{\varepsilon_s}\frac{\partial}{\partial y_i^s},\;\;\mathbf{x} \in\Omega_s
\end{equation}
which will be extensively used in the sequel. Hence, $k_{ij}^{\bm{\varepsilon }}({\mathbf{x}})$
can be changed into $k_{ij}({\mathbf{y}}^s)$ for $\mathbf{x}\in \Omega_s$. Being similar to \cite{R3,R4,R5,R6}, we make the following assumptions:
\begin{enumerate}[(A)]
\item[(A)]
$k_{ij}({\mathbf{y}}^s)$ is a scalar function belonging to $L^\infty (\Omega)$ and function $k_{ij}({\mathbf{y}}^s)$ is $1-$periodic (with $Q^s=(0,1)^{\mathcal{N}}$) for $\mathbf{x}\in \Omega_s$.
\item[(B)]
functions $k_{ij}({\mathbf{y}}^s)$ is symmetric and there exist two constants $0<\gamma_0\leq\gamma_1$ such that
\begin{displaymath}
\gamma_0|\bm{\xi}|^2\leq k_{ij}^{\bm{\varepsilon }}({\mathbf{x}})\xi_i\xi_j  \le\gamma_1|\bm{\xi}|^2,\;\;\mathbf{x}\in\Omega
\end{displaymath}
for all vectors $\bm{\xi}=(\xi_i)\in \mathbb{R}^\mathcal{N}$ and all symmetric matrix $\{\eta_{ij}\}\in \mathbb{R}^\mathcal{N}\times \mathbb{R}^\mathcal{N}$.
\item[(C)]
$h(\mathbf{x})\in L^2(\Omega),\widehat{T}(\mathbf{x})\in L^2(\Omega),\bar{q}(\mathbf{x})\in L^2(\Omega)$.
\end{enumerate}

Now, we give the specific construction process of the FOTS solution and SOTS solution for problem (1). To the problem (1), we assume that $T^{\bm{\varepsilon}}(\mathbf{x})$ can be formally expanded as follows:
\begin{equation}
T^{\bm{\varepsilon}}(\mathbf{x}) = {T^{(0)}}(\mathbf{x},\mathbf{y}^s) + \varepsilon_s {T^{(1)}}(\mathbf{x},\mathbf{y}^s) + {\varepsilon_s^2}{T^{(2)}}(\mathbf{x},\mathbf{y}^s) + {\rm O}({\varepsilon_s^3}),\;\;\mathbf{x} \in\Omega_s
\end{equation}
Substituting (3) into (1) and by virtue of the chain rule (2), we have
\begin{equation}
\begin{aligned}
&- \varepsilon _s^{ - 2}\frac{\partial }{\partial y_i^s}\Big( {{k_{ij}}({\bf{y}}^s)\frac{\partial T^{(0)}}{\partial y_j^s}} \Big) - \varepsilon_s^{ - 1}\frac{\partial }{\partial y_i^s}\Big( {k_{ij}({\bf{y}}^s)\big( {\frac{\partial T^{(0)}}{\partial {x_j}} + \frac{\partial T^{(1)}}{\partial y_j^s}} \big)} \Big)\\
&- \varepsilon _s^{ - 1}\frac{\partial }{\partial {x_i}}\Big( {k_{ij}({\bf{y}}^s)\frac{\partial {T^{(0)}}}{\partial y_j^s}} \Big)- \varepsilon _s^0\frac{\partial }{\partial y_i^s}\Big( {k_{ij}({\bf{y}}^s)\big( {\frac{\partial {T^{(1)}}}{\partial {x_j}} + \frac{\partial {T^{(2)}}}{\partial y_j^s}} \big)} \Big)\\
&- \varepsilon _s^0\frac{\partial }{\partial {x_i}}\Big( {k_{ij}({\bf{y}}^s)\big( {\frac{\partial {T^{(0)}}}{\partial {x_j}} + \frac{\partial {T^{(1)}}}{\partial y_j^s}} \big)} \Big)+ O({\varepsilon _s})= h,\;\;\mathbf{x} \in\Omega_s
\end{aligned}
\end{equation}
From (4), a series of equations in $\Omega_s$ are derived by matching terms of the same order of $\varepsilon_s$ according to the classical procedure of AHM \cite{R11,R12}
\begin{equation}
O(\varepsilon _s^{ - 2}):-\frac{\partial }{\partial y_i^s}\Big( {k_{ij}({\bf{y}}^s)\frac{\partial {T^{(0)}}}{\partial y_j^s}} \Big) = 0
\end{equation}
\begin{equation}
O(\varepsilon _s^{ - 1}):-\frac{\partial }{\partial y_i^s}\Big( {k_{ij}({\bf{y}}^s)\big( {\frac{\partial {T^{(0)}}}{\partial {x_j}} + \frac{\partial {T^{(1)}}}{\partial y_j^s}} \big)} \Big) - \varepsilon _s^{ - 1}\frac{\partial }{\partial {x_i}}\Big( {k_{ij}({\bf{y}}^s)\frac{\partial {T^{(0)}}}{\partial y_j^s}} \Big) = 0
\end{equation}
\begin{equation}
\begin{aligned}
O(\varepsilon _s^0):- \frac{\partial }{{\partial {x_i}}}\Big( {{k_{ij}}({{\bf{y}}^s})\big( {\frac{{\partial {T^{(0)}}}}{{\partial {x_j}}} + \frac{{\partial {T^{(1)}}}}{{\partial y_j^s}}} \big)} \Big)- \frac{\partial }{{\partial y_i^s}}\Big( {{k_{ij}}({{\bf{y}}^s})\big( {\frac{{\partial {T^{(1)}}}}{{\partial {x_j}}} + \frac{{\partial {T^{(2)}}}}{{\partial y_j^s}}} \big)} \Big) = h
\end{aligned}
\end{equation}
From (5) we can acquire that $T^{(0)}(\mathbf{x},\mathbf{y}^s)$ are independent of the micro-scale variable $\mathbf{y^s}$, namely
\begin{equation}
T^{(0)}(\mathbf{x},\mathbf{y}^s) =T^{(0)}(\mathbf{x}),\;\;\mathbf{x} \in\Omega_s
\end{equation}
Subsequently, (6) can be further simplified as the following equations
\begin{equation}
\frac{\partial }{\partial y_i^s}\Big( {k_{ij}({\bf{y}}^s)\frac{\partial {T^{(1)}}}{\partial y_j^s}}\Big)=-\frac{\partial }{\partial y_i^s}\Big( {k_{ij}({\bf{y}}^s)\frac{\partial T^{(0)}}{\partial {x_j}}} \Big),\;\;\mathbf{x} \in\Omega_s
\end{equation}
According to (9), we construct
\begin{equation}
T^{(1)}({\bf{x}},{{\bf{y}}^s}) = M_{{\alpha _1}}({\bf{y}}^s)\frac{\partial {T^{(0)}}}{\partial {x_{{\alpha _1}}}},\;\;\mathbf{x} \in\Omega_s
\end{equation}
where $M_{{\alpha _1}}({\bf{y}}^s)$ are the first-order auxiliary cell functions defined in unit cell $Q^s$.

Now, substituting (10) into (9), the following equations with homogeneous Dirichlet boundary condition are obtained after simplification and calculation
\begin{equation}
\left\{ \begin{aligned}
&\frac{\partial }{\partial y_i^s}\Big( k_{ij}({\bf{y}}^s)\frac{\partial {M_{\alpha _1}({\bf{y}}_s})}{{\partial y_j^s}} \Big) =-\frac{\partial k_{i\alpha_1}({\bf{y}}^s)}{\partial y_i^s},&\;\;&{\bf{y}}^s \in {Q^s}\\
&M_{{\alpha _1}}({\bf{y}}^s) = 0,&\;\;&{\bf{y}}^s \in \partial {Q^s}
\end{aligned} \right.
\end{equation}
After that, we make the volume integral to both sides of (7) on the unit cell $Q^s$ and using the Gauss theorem on (7)
\begin{equation}
\left\{ \begin{aligned}
&- \frac{\partial }{{\partial {x_i}}}\left( {\widehat k_{ij}^s\frac{{\partial {T^{(0)}}({\bf{x}})}}{{\partial {x_j}}}} \right)= h(\mathbf{x}),\;\;\text{in}\;\;\Omega=\bigcup\limits_{s = 1}^K{\Omega _s},\;\;\bm{\varepsilon}=(\varepsilon _1,\cdots\varepsilon_K),\\
&T^{(0)}(\mathbf{x}) = \widehat T(\mathbf{x})\;\;\;\text{on}\;\;\partial {\Omega_T},\\
&{\widehat k_{ij}^s\frac{\partial T^{(0)}(\mathbf{x})}{\partial {x_j}}}{n_i} = \bar q(\mathbf{x})\;\;\;\text{on}\;\;\partial {\Omega_q}
\end{aligned} \right.
\end{equation}
where the homogenized material parameters are defined as follows
\begin{equation}
{\widehat k_{ij}^s = \frac{1}{{\left| {{Q^s}} \right|}}\int_{{Q^s}} {\left( {{k_{ij}}({{\bf{y}}^s}) + {k_{i{\alpha _1}}}({{\bf{y}}^s})\frac{{\partial {M_j}({{\bf{y}}^s})}}{{\partial y_{{\alpha _1}}^s}}} \right)d{Q^s}} },\;\;\mathbf{x} \in\Omega_s
\end{equation}
Now, we start to solve the vital second-order auxiliary cell functions. Firstly, the following equations are obtained by subtracting (7) from (12)
\begin{equation}
\begin{aligned}
&\frac{\partial}{\partial y_i^s}\Big( {k_{ij}({\bf{y}}^s){\frac{\partial T^{(2)}}{\partial y_j^s}}} \Big)=\Big[ { {\widehat k}_{\alpha_1\alpha_2}^s - k_{\alpha_1\alpha_2}({\bf{y}}^s)}-{ \frac{\partial}{\partial y_i^s}\big( {k_{i\alpha_2}({\bf{y}}^s){M_{\alpha_1}({\bf{y}}^s)}} \big)}\\
&{-{k_{\alpha_1j}({\bf{y}}^s)\frac{\partial M_{\alpha_2}({\bf{y}}^s)}{\partial y_j^s}}} \Big]\frac{\partial^2 T^{(0)}}{\partial x_{\alpha_1}\partial x_{\alpha_2}},\;\;\mathbf{x} \in\Omega_s
\end{aligned}
\end{equation}
According to (14), the specific form of ${T^{(2)}}({\bf{x}},{{\bf{y}}^s})$ is constructed as follows
\begin{equation}
{T^{(2)}}({\bf{x}},{{\bf{y}}^s}) = {M_{{\alpha _1}{\alpha _2}}}({{\bf{y}}^s})\frac{{{\partial ^2}{T^{(0)}}}}{{\partial {x_{{\alpha _1}}}\partial {x_{{\alpha _2}}}}},\;\;\mathbf{x} \in\Omega_s
\end{equation}
where $M_{\alpha_1\alpha_2}({\bf{y}}^s)$ are the second-order auxiliary cell functions defined in unit cell $Q^s$.

Substituting (15) into (14), the following equations, which are attached with the homogeneous Dirichlet boundary condition, are derived as follows
\begin{equation}
\left\{
\begin{aligned}
&\frac{\partial}{\partial y_i^s}\Big( { k_{ij}({\bf{y}}^s){\frac{\partial {M_{\alpha_1\alpha_2}({\bf{y}}^s)}}{\partial y_j^s}}} \Big) = { {\widehat k}_{\alpha_1\alpha_2}^s - k_{\alpha_1\alpha_2}({\bf{y}}^s)}\\
&-{ \frac{\partial}{\partial y_i^s}\big( {k_{i\alpha_2}({\bf{y}}^s){M_{\alpha_1}({\bf{y}}^s)}} \big)-{k_{\alpha_1j}({\bf{y}}^s)\frac{\partial M_{\alpha_2}({\bf{y}}^s)}{\partial y_j^s}}},&\;\;&{\bf{y}}^s \in {Q^s}\\
&{M_{\alpha_1\alpha_2}}({\bf{y}}^s)=0,&\;\;&{\bf{y}}^s \in\partial {Q^s}
\end{aligned} \right.
\end{equation}
In a summary, we can get the following theorem for multiscale problem (1).
\begin{thm}
The heat conduction problems of composite structures with diverse periodic configurations in different subdomains have SOTS asymptotic expansion solutions as follows
\begin{equation}
{\begin{aligned}
{T^{\bm{\varepsilon}} }(\mathbf{x})&\cong T^{(0)}+\varepsilon_s{M_{\alpha_1}}({\bf{y}}^s)\frac{\partial T^{(0)}}{\partial x_{\alpha_1}}+ \varepsilon^2_s{M_{\alpha_1\alpha_2}}({\bf{y}}^s)\frac{\partial^2 T^{(0)}}{\partial x_{\alpha_1}\partial x_{\alpha_2}},\\
&\;\;\mathbf{x} \in\Omega_s,\;\;s=1,\cdots,K
\end{aligned}}
\end{equation}
\end{thm}
\noindent where $T^{(0)}$ is the solution of the homogenized problem (12). $M_{\alpha_1}$ is the first-order auxiliary cell functions defined by (11). $M_{\alpha_1\alpha_2}$ is the second-order auxiliary cell functions defined by (16).
\section{Error analysis of multiscale approximate solutions}
In this section, the detailed error analysis of FOTS solutions and SOTS solutions in the pointwise sense is given. Firstly, we denote the FOTS solutions $T^{(1\bm{\varepsilon})}$ and SOTS solutions $T^{(2\bm{\varepsilon})}$ for governing equation as follows:
\begin{equation}
T^{(1\bm{\varepsilon})}=T^{(0)}+\varepsilon_s T^{(1)},\;T^{(2\bm{\varepsilon})}=T^{(0)}+\varepsilon_s T^{(1)}+\varepsilon^2_s T^{(2)},\;\;\mathbf{x} \in\Omega_s
\end{equation}
Then, define the following residual functions for the FOTS solutions and SOTS solutions
\begin{equation}
T^{(1\bm{\varepsilon})}_{\Delta}=T^{\bm{\varepsilon}}-T^{(1\bm{\varepsilon})},\;T^{(2\bm{\varepsilon})}_{\Delta}=T^{\bm{\varepsilon}}-T^{(2\bm{\varepsilon})}
\end{equation}
Before giving the detailed analysis procedure, we need to make some assumptions about multiscale problem (1). Assume that $\Omega\subset\mathbb{R}^\mathcal{N}$ is a bounded domain and each $\Omega_s$ is the union of entire periodic cells, i.e. $\bar{\Omega}_s=\cup_{\mathbf{z}\in T_{\varepsilon_s}}\varepsilon(\mathbf{z}+\bar{Q}^s)$, where the index set $T_{\varepsilon_s}=\{\mathbf{z}=(z_1,\cdots,z_\mathcal{N})\in Z^\mathcal{N},\varepsilon(\mathbf{z}+\bar{Q}^s)\subset \bar{\Omega}_s\}$. Besides, let $E_\mathbf{z}^s=\varepsilon(\mathbf{z}+Q^s)$ and $\partial E_\mathbf{z}^s$ be the boundary of $E_\mathbf{z}$.
\subsection{Error analysis in the pointwise sense}
To compare $T^{(1\bm{\varepsilon})}$ with the original solution $T^{\bm{\varepsilon}}$, substituting the residual function $T^{(1\bm{\varepsilon})}_{\Delta}$ into (1) and using (11) and (12), the following residual equation of FOTS solutions is obtained which holds in the distribution sense:
\begin{equation}
\left\{ \begin{aligned}
&- \frac{\partial }{{\partial {x_i}}}\Big( {k_{ij}^{\bm{\varepsilon }}({\mathbf{x}})\frac{{\partial {T_\Delta^{(1\bm{\varepsilon })} }}}{{\partial {x_j}}}}\Big)=F_{0}(\mathbf{x},\mathbf{y}^s)\\
&\quad\quad\quad+\varepsilon_s F_{1}(\mathbf{x},\mathbf{y}^s),\;\;\mathbf{x} \in\Omega_s,\;\;s=1,\cdots,K,\\
&T_{\Delta}^{(1\bm{\varepsilon })}(\mathbf{x})=0\;\;\;\text{on}\;\;\partial\Omega_T,\\
&{k_{ij}^{\bm{\varepsilon }}({\mathbf{x}})\frac{{\partial {T_\Delta^{(1\bm{\varepsilon })} }}}{{\partial {x_j}}}}n_i=\Re_{1i}n_i\;\;\;\text{on}\;\;\partial {\Omega_q}
\end{aligned} \right.
\end{equation}
where the detailed forms of $F_{0}(\mathbf{x},\mathbf{y}^s)$ and $F_{1}(\mathbf{x},\mathbf{y}^s)$ are listed as follows:
\begin{equation}
\begin{aligned}
F_{0}(\mathbf{x},\mathbf{y}^s)&=-\frac{\partial }{{\partial {x_i}}}\big( {{\widehat k_{ij}^s}\frac{{\partial {T^{(0)}}}}{{\partial {x_j}}}}\big)+\frac{\partial }{{\partial {x_i}}}\big( {{ k_{ij}}\frac{{\partial {T^{(0)}}}}{{\partial {x_j}}}}\big)\\
&+\frac{\partial }{{\partial {x_i}}}\Big[ {{ k_{ij}}\frac{{\partial }}{{\partial {y_j^s}}}}\big(M_{\alpha_1}\frac{\partial T^{(0)}}{\partial x_{\alpha_1}}\big)\Big]
+\frac{\partial }{{\partial {y_i^s}}}\Big[ {{ k_{ij}}\frac{{\partial }}{{\partial {x_j}}}}\big(M_{\alpha_1}\frac{\partial T^{(0)}}{\partial x_{\alpha_1}}\big)\Big]
\end{aligned}
\end{equation}
\begin{equation}
F_{1}(\mathbf{x},\mathbf{y}^s)=\frac{\partial }{{\partial {x_i}}}\Big[ {{ k_{ij}}\frac{{\partial }}{{\partial {x_j}}}}\big(M_{\alpha_1}\frac{\partial T^{(0)}}{\partial x_{\alpha_1}}\big)\Big]
\end{equation}
Then substituting $T^{(2\bm{\varepsilon })}_{\Delta}$ into (1) and by virtue of (11), (12) and (16), we obtain the following residual equation of the SOTS solutions which holds in the distribution sense:
\begin{equation}
\left\{ \begin{aligned}
&- \frac{\partial }{{\partial {x_i}}}\Big( {k_{ij}^{\bm{\varepsilon }}({\mathbf{x}})\frac{{\partial {T_\Delta^{(2\bm{\varepsilon })} }}}{{\partial {x_j}}}}\Big)=\varepsilon_s G(\mathbf{x},\mathbf{y}^s),\;\;\mathbf{x} \in\Omega_s,\;\;s=1,\cdots,K,\\
&T_{\Delta}^{(2\bm{\varepsilon })}(\mathbf{x})=0\;\;\;\text{on}\;\;\partial\Omega_T,\\
&{k_{ij}^{\bm{\varepsilon }}({\mathbf{x}})\frac{{\partial {T_\Delta^{(2\bm{\varepsilon })} }}}{{\partial {x_j}}}}n_i=\Re_{2i}n_i\;\;\;\text{on}\;\;\partial {\Omega_q}
\end{aligned} \right.
\end{equation}
where the detailed form of $G(\mathbf{x},\mathbf{y}^s)$ are given as follows:
\begin{equation}
\begin{aligned}
G(\mathbf{x},\mathbf{y}^s)&=\frac{\partial}{\partial y_i^s}\Big[k_{ij}\frac{\partial}{\partial x_j}\big( {M_{\alpha_1\alpha_2}}\frac{\partial^2 T^{(0)}}{\partial x_{\alpha_1}\partial x_{\alpha_2}}\big) \Big]+\frac{\partial}{\partial x_i}\Big[k_{ij}\frac{\partial}{\partial y_j^s}\big( {M_{\alpha_1\alpha_2}}\frac{\partial^2 T^{(0)}}{\partial x_{\alpha_1}\partial x_{\alpha_2}}\big)\Big]\\
&+\varepsilon\frac{\partial}{\partial x_i}\Big[k_{ij}\frac{\partial}{\partial x_j}\big( R_{\alpha_1}\frac{\partial T^{(0)}}{\partial x_{\alpha_1}}+{M_{\alpha_1\alpha_2}}\frac{\partial^2 T^{(0)}}{\partial x_{\alpha_1}\partial x_{\alpha_2}}\big) \Big]
\end{aligned}
\end{equation}
Now we can give a conclusion about the error analysis in the pointwise sense. From the residual equation (20), one can easily see that the residual of FOTS solutions is order $O(1)$ in the pointwise sense due to the term $F_0(\mathbf{x},\mathbf{y}^s)$. In addition, it is clear to see that the residual of SOTS solutions is order $O(\varepsilon)$ in the pointwise sense from the residual equation (23). This means that SOTS solutions can satisfy the original equation (1) in the pointwise sense. Thus even $\varepsilon$ is a small constant, the SOTS solutions can still provide the required accuracy of engineering calculation and capture the micro-scale oscillating behavior of composite materials. This is the main reason and motivation to develop the SOTS solutions.
\subsection{Main convergence theorem and its proof}
The error estimates for the SOTS approximate solutions in the integral sense is presented in this subsection. It is known to all that the classical auxiliary cell functions are defined with periodic boundary conditions and have enough regularity on the boundary of unit cell $Q^s$ \cite{R12,R13,R14,R15}. In the case of auxiliary cell functions defined with homogeneous Dirichlet boundary condition, their normal derivatives are only continuous on the boundary of $Q^s$ under geometric symmetry and regularity assumptions on material property parameters. So we firstly give some hypotheses similar to literature as follows \cite{R13,R14,R15}:
\begin{enumerate}[(i)]
\item
$k_{ij}(\mathbf{y}^s)$ is a function with piecewise constants in all $Q^s$.
\item
Let $\Delta_1\ldots\Delta_\mathcal{N}(\mathcal{N}=2,3)$ be the middle hyperplanes of the reference cell $Q^s=(0,1)^{\mathcal{N}}$. Assume that $k_{ii}(\mathbf{y}^s)$ is symmetric with respect to $\Delta_1\ldots\Delta_\mathcal{N}$ and $k_{ij}(\mathbf{y}^s)$ is anti-symmetric with respect to $\Delta_1\ldots\Delta_\mathcal{N}$ in $Q^s$ for $\mathbf{x}\in \Omega_s$.
\end{enumerate}
\begin{lem}
Denote operator $\displaystyle\sigma_{Q^s}=n_i k_{ij}(\mathbf{y}^s)\frac{\partial}{\partial y_j^s}$ for $\mathbf{x}\in \Omega_s$. Then under assumptions (A)-(B) and (i)-(ii), the normal derivatives $\sigma_{Q^s}(M_{\alpha_1})$ and $\sigma_{Q^s}(M_{\alpha_1\alpha_2})$ can be proved to be continuous on the boundary of unit cell $Q^s$ by using the same method in Refs. \cite{R13,R14,R15}.
\end{lem}
\begin{thm}
Assume that $\Omega\subset\mathbb{R}^\mathcal{N}$ is a bounded domain and each $\Omega_s$ is the union of entire periodic cells, i.e. $\bar{\Omega}_s=\cup_{\mathbf{z}\in T_{\varepsilon_s}}\varepsilon(\mathbf{z}+\bar{Q}^s)$, where the index set $T_{\varepsilon_s}=\{\mathbf{z}=(z_1,\cdots,z_\mathcal{N})\in Z^\mathcal{N},\varepsilon(\mathbf{z}+\bar{Q}^s)\subset \bar{\Omega}_s\}$. Let $T^{\bm{\varepsilon}}(\mathbf{x})$ be the weak solution of multiscale problem (1), $T^{(0)}(\mathbf{x})$ is the solution of associated homogenized problem (12). $T^{(2\bm{\varepsilon})}(\mathbf{x})$ is the SOTS approximate solution stated in Theorem 1. Under the aforementioned assumptions (A)-(C), (i)-(ii), and Lemma 1, we obtain the following error estimate
\begin{equation}
\big\|T^{\bm{\varepsilon}}(\mathbf{x})-T^{(2\bm{\varepsilon})}(\mathbf{x})\big\|_{H^1(\Omega)}\leq C\varepsilon_{max}^{\frac{1}{2}},\;\;\varepsilon_{max}=max\{\varepsilon _1,\cdots\varepsilon_K\}
\end{equation}
\end{thm}
where $C$ is a positive constant independent of $\bm{\varepsilon}$, but dependent of $\Omega$.\\
$\mathbf{Proof:}$ Firstly, the following equality can be obtained from (2) and (17)
\begin{equation}
\begin{aligned}
\sigma_{T}(T^{(2\bm{\varepsilon})})=&n_j k_{ij}^{\bm{\varepsilon }}({\mathbf{x}})\frac{\partial T^{(2\bm{\varepsilon})} }{\partial x_j}\\
=&n_j k_{ij}(\mathbf{y}^s)\big(\frac{\partial}{\partial x_j}+\frac{1}{\varepsilon_s}\frac{\partial}{\partial y_j^s}\big)\Big[T^{(0)}+\varepsilon_s{M_{\alpha_1}}\frac{\partial T^{(0)}}{\partial x_{\alpha_1}}+ \varepsilon_s^2{M_{\alpha_1\alpha_2}}\frac{\partial^2 T^{(0)}}{\partial x_{\alpha_1}\partial x_{\alpha_2}}\Big]\\
=&n_j k_{ij}(\mathbf{y}^s)\frac{\partial}{\partial x_j}\Big[T^{(0)}+\varepsilon_s{M_{\alpha_1}}\frac{\partial T^{(0)}}{\partial x_{\alpha_1}}+ \varepsilon_s^2{M_{\alpha_1\alpha_2}}\frac{\partial^2 T^{(0)}}{\partial x_{\alpha_1}\partial x_{\alpha_2}}\Big]\\
&+\Big[{\sigma_{Q^s}(M_{\alpha_1})}\frac{\partial T^{(0)}}{\partial x_{\alpha_1}}+ \varepsilon_s{\sigma_{Q^s}(M_{\alpha_1\alpha_2})}\frac{\partial^2 T^{(0)}}{\partial x_{\alpha_1}\partial x_{\alpha_2}}\Big],\;\;\mathbf{x} \in\Omega_s
\end{aligned}
\end{equation}
Secondly, we use the residual equation (23) to complete the error estimate. Multiplying by $\displaystyle T_{\Delta}^{(2\bm{\varepsilon})}$ on both sides of (23) and integrating on each $\Omega_s$, then the following equations are derived by summing up of all $\Omega_s$
\begin{equation}
-\sum_{s=1}^{K}\int_{\Omega_s} \frac{\partial}{\partial x_i}\Big( {{k_{ij}}(\mathbf{y}^s)\frac{{\partial {T_\Delta^{(2\bm{\varepsilon})} }}}{\partial {x_j}}}\Big)T_\Delta^{(2\bm{\varepsilon})}d\Omega_s
=\sum_{s=1}^{K}\int_{\Omega_s}\varepsilon_s G(\mathbf{x},\mathbf{y}^s)T_\Delta^{(2\bm{\varepsilon})}d\Omega_s
\end{equation}
Using Green's formula and integrating by parts on (27), (27) can be simplified as follows
\begin{equation}
\begin{aligned}
&\sum_{s=1}^{K}\int_{\Omega_s} {{k_{ij}}(\mathbf{y}^s)\frac{{\partial {T_\Delta^{(2\bm{\varepsilon})} }}}{\partial {x_j}}}\frac{{\partial {T_\Delta^{(2\bm{\varepsilon})} }}}{\partial {x_i}}d\Omega_s
=\sum_{s=1}^{K}\int_{\Omega_s}\varepsilon_s G(\mathbf{x},\mathbf{y}^s)T_\Delta^{(2\bm{\varepsilon})}d\Omega_s\\
&+\int_{\partial\Omega_q}\Re_{2i}n_iT_\Delta^{(2\bm{\varepsilon})}ds+\sum_{s=1}^{K}{\int_{\cup_{\mathbf{z}\in T_{\varepsilon_s}}\partial E_\mathbf{z}^s}}{\varphi}T_{\Delta}^{(2\bm{\varepsilon})}d\Gamma_\mathbf{y}
\end{aligned}
\end{equation}
where $\varphi$ results from using the Green's formula on $\partial E_z^s$.

Combining (26) and Lemma 1, we can obtain
\begin{equation}
\begin{aligned}
\big \langle {\varphi},T_{\Delta}^{(2\bm{\varepsilon})}\big \rangle= &\sum_{s=1}^{K}{\int_{\cup_{\mathbf{z}\in T_{\varepsilon_s}}\partial E_\mathbf{z}^s}}{\varphi}T_{\Delta}^{(2\bm{\varepsilon})}d\Gamma_\mathbf{y}\\
=&\sum_{s=1}^{K}\sum\limits _{\mathbf{z}\in T_{\varepsilon_s}}\int_{\partial E_\mathbf{z}^s}\sigma_T(T^{\bm{\varepsilon}}-T^{(2\bm{\varepsilon})})T_{\Delta}^{(2\bm{\varepsilon})}d\Gamma_\mathbf{y}\\
=&-\sum_{s=1}^{K}\sum\limits _{\mathbf{z}\in T_{\varepsilon_s}}\int_{\partial E_\mathbf{z}^s}\sigma_T(T^{(2\bm{\varepsilon})})T_{\Delta}^{(2\bm{\varepsilon})}d\Gamma_\mathbf{y}=0
\end{aligned}
\end{equation}
Afterwards, it is easy to derive the following identity by substituting (29) into (28)
\begin{equation}
\begin{aligned}
&\sum_{s=1}^{K}\int_{\Omega_s} {{k_{ij}}(\mathbf{y}^s)\frac{{\partial {T_\Delta^{(2\bm{\varepsilon})} }}}{\partial {x_j}}}\frac{{\partial {T_\Delta^{(2\bm{\varepsilon})} }}}{\partial {x_i}}d\Omega_s
=\sum_{s=1}^{K}\int_{\Omega_s}\varepsilon_s G(\mathbf{x},\mathbf{y}^s)T_\Delta^{(2\bm{\varepsilon})}d\Omega_s\\
&+\int_{\partial\Omega_q}\Re_{2i}n_iT_\Delta^{(2\bm{\varepsilon})}ds
\end{aligned}
\end{equation}
We underline that the equality (30) is vital to obtain the error estimate (25).

Now, applying the Poincar$\acute{e}$-Friedrichs inequality to the left side of (30), we gets
\begin{equation}
\Big|\sum_{s=1}^{K}\int_{\Omega_s} {{k_{ij}}(\mathbf{y}^s)\frac{{\partial {T_\Delta^{(2\bm{\varepsilon})} }}}{\partial {x_j}}}\frac{{\partial {T_\Delta^{(2\bm{\varepsilon})} }}}{\partial {x_i}}d\Omega_s\Big| \geq C\left \|T_\Delta^{(2\bm{\varepsilon})}\right\|_{H^1(\Omega)}^2
\end{equation}
After that, making use of the Schwarz's inequality, theorem 1.2 and lemma 2.2 in \cite{R11}, the following inequality is obtained by transforming the right side of (30):
\begin{equation}
\begin{aligned}
&\Big|\sum_{s=1}^{K}\int_{\Omega_s}\varepsilon_s G(\mathbf{x},\mathbf{y}^s)T_\Delta^{(2\bm{\varepsilon})}d\Omega_s
+\int_{\partial\Omega_q}\Re_{2i}n_iT_\Delta^{(2\bm{\varepsilon})}ds\Big|\\
&\leq\Big|\sum_{s=1}^{K}\int_{\Omega_s}\varepsilon_s G(\mathbf{x},\mathbf{y}^s)T_\Delta^{(2\bm{\varepsilon})}d\Omega_s\Big|+\Big|\int_{\partial\Omega_q}\Re_{2i}n_iT_\Delta^{(2\bm{\varepsilon})}ds\Big|\\
&\leq \sum_{s=1}^{K}\left\|\varepsilon_s G(\mathbf{x},\mathbf{y}^s)\right\|_{L^2(\Omega_s)} \left\|T_\Delta^{(2\bm{\varepsilon})}\right\|_{L^2(\Omega_s)}+C\varepsilon^{\frac{1}{2}}_{max}\left\|T_\Delta^{(2\bm{\varepsilon})}\right\|_{H^1(\Omega)}\\
&\leq C\varepsilon_{max}\left\|T_\Delta^{(2\bm{\varepsilon})}\right\|_{H^1(\Omega)}+C\varepsilon^{\frac{1}{2}}_{max}\left\|T_\Delta^{(2\bm{\varepsilon})}\right\|_{H^1(\Omega)}
\end{aligned}
\end{equation}
Combining (31) and (32) together, it follows that
\begin{equation}
\left\|T_\Delta^{(2\bm{\varepsilon})}\right\|_{H^1(\Omega)}\leq C\varepsilon_{max}+C\varepsilon^{\frac{1}{2}}_{max}
\end{equation}
Finally, it is obvious that we verify
\begin{equation}
\big\| T^{\bm{\varepsilon}}(\mathbf{x})-T^{(2\bm{\varepsilon})}(\mathbf{x})\big\|_{H^1(\Omega)}=\left\|T^{(2\bm{\varepsilon})}_{\Delta}\right\|_{H^1(\Omega)}\leq C\varepsilon^{\frac{1}{2}}_{max}.
\end{equation}
where $C$ denotes a positive generic constant and has different values in different places in this paper.
\section{Second-order two-scale numerical algorithm}
In this section, we give the detailed SOTS numerical algorithm for the multiscale problem (1). From the SOTS analysis of multiscale problem (1), we underline that the auxiliary cell problems (11) and (16) have different solutions on different unit cell $Q^s$. It is totally different from classical AHM. Based on the above-mentioned analysis, we present the following SOTS numerical algorithm for model problem (1), which is based on FDM in time direction and FEM in spatial region. The detailed algorithm procedures are listed as follows:
\begin{enumerate}[(1)]
\item
Define the geometric structure of the unit cell $Q^s=(0,1)^\mathcal{N}$ and homogenized macroscopic region $\Omega$ in $\mathbb{R}^\mathcal{N}$, and verify the material parameters of composite materials. Then, generate the triangular finite element mesh in $\mathbb{R}^2$ or tetrahedral mesh in $\mathbb{R}^3$. Let $J^{h_1}=\{K\}$ and $J^{h_0}=\{e\}$ be a regular family of triangles or tetrahedra of the unit cell $Q^s$ and the homogenized macroscopic region $\Omega$, respectively, where $h_1=$max$_K\{h_K\}$ and $h_0=$max$_e\{h_e\}$. And define the linear conforming finite element spaces $V_{h_1}(Q^s)=\{\nu\in C^0(\bar{Q}^s):\nu\mid_{\partial Q^s}=0,\nu\mid_{K}\in P_1(K)\}\subset H_0^1(Q^s)$ and $V_{h_0}(\Omega)=\{\nu\in C^0(\bar{\Omega}):\nu\mid_{\partial\Omega_{T}}=0,\nu\mid_{e}\in P_1(e)\}\subset H^1(\Omega)$ for the above two regions, respectively.
\item
Solve the first-order auxiliary cell problems (11) on $V_{h_1}(Q^s)$ corresponding to different unit cell $Q^s$. And the homogenized material parameters ${\widehat S}^s$ and $\widehat k_{ij}^s$ are evaluated by making integral of (13) corresponding to different unit cell $Q^s$. After that, the homogenized material parameters on each nodes of $V_{h_0}(\Omega)$ can be determined by identifying the subdomain of their coordinates.
\item
Then, in turn, solve the homogenized equations (12) in the macroscopic $\Omega$ to obtain the homogenized temperature $T^{(0)}$. The following FEM scheme is adopted to compute the homogenized heat conduction problem (12)
\begin{equation}
\left\{\begin{aligned}
&\sum_{s=1}^{K}\int_{\Omega}\widehat{k}_{ij}^s\frac{\partial T^{(0)}}{\partial x_j}\frac{\partial \widetilde{\varphi}^{h_0}}{\partial x_i}d\Omega_s=\int_{\Omega} h\widetilde{\varphi}^{h_0}d\Omega\\
&+\int_{\partial\Omega_q}\bar{q}(\mathbf{x})\widetilde{\varphi}^{h_0}ds,\;\forall\widetilde{\varphi}^{h_0}\in V_{h_0}(\Omega),\\
&T^{(0)}({\mathbf{x}})=\widehat{T}({\mathbf{x}})\;\;\;\text{on}\;\;\partial\Omega_{T}
\end{aligned}\right.
\end{equation}
\item
Using the same mesh as first-order auxiliary cell problems, the second-order auxiliary cell problems (16), which correspond to different unit cell $Q^s$, are solved on $V_{h_1}(Q^s)$, respectively.
\item
For arbitrary point $\mathbf{x}\in \Omega$, we use the interpolation method to get the corresponding values of first-order auxiliary cell functions, second-order auxiliary cell functions and homogenized solutions. The spatial derivatives $\displaystyle\frac{\partial T^{(0)}}{\partial x_{\alpha_1}}$ and $\displaystyle\frac{\partial^2 T^{(0)}}{\partial x_{\alpha_1}\partial x_{\alpha_2}}$ are evaluated by the average technique on relative elements \cite{R9,R15}. Then, the displacement field temperature field $T^{(2\bm{\varepsilon})}(\mathbf{x})$ can be solved by the formula (17). Moreover, we can still use the higher-order interpolation method and post-processing technique to get the high-precision SOTS solutions \cite{R15,R16}.
\end{enumerate}
\section{Numerical examples}
In this section, two numerical examples are given to check the validity and feasibility of the SOTS numerical algorithm we developed. Since it is difficult to find the analytical solutions for the multiscale problem (1), we replace $T^{\bm{\varepsilon}}(\mathbf{x})$ with $T_e(\mathbf{x})$ which is precise FEM solutions for multiscale problem (1) on a very fine mesh. Without confusion, some notations are introduced as follows:
\begin{equation}
Terror0=\frac{||T_e-T^{(0)}||_{L^2}}{||T_e||_{L^2}},
Terror1=\frac{||T_e-T^{(1\bm{\varepsilon})}||_{L^2}}{||T_e||_{L^2}},
Terror2=\frac{||T_e-T^{(2\bm{\varepsilon})}||_{L^2}}{||T_e||_{L^2}}
\end{equation}
\begin{equation}
TError0=\frac{|T_e-T^{(0)}|_{H^1}}{|T_e|_{H^1}},
TError1=\frac{|T_e-T^{(1\bm{\varepsilon})}|_{H^1}}{|T_e|_{H^1}},
TError2=\frac{|T_e-T^{(2\bm{\varepsilon})}|_{H^1}}{|T_e|_{H^1}}
\end{equation}
\subsection{Example 1: different configurations with same material}
In this example, a 2D composite structures with two basic periodic configurations in four subdomains is considered. The macrostructure $\Omega$ is shown in Fig. 1, where $\Omega=(x,y)=[0,2]\times[0,2]cm^2$. Assume that $\Omega=\Omega_1\cup\Omega_2\cup\Omega_3\cup\Omega_4$ where $\Omega_1=[0,1]\times[0,1]cm^2$, $\Omega_2=[1,2]\times[0,1]cm^2$, $\Omega_3=[1,2]\times[1,2]cm^2$ and $\Omega_4=[0,1]\times[1,2]cm^2$. Moreover, $\Omega_1$ and $\Omega_3$ have the same unit cell $Q^1$ and periodic unit cell size $\displaystyle\varepsilon_1=\frac{1}{6}$. $\Omega_2$ and $\Omega_4$ have the same unit cell $Q^2$ and periodic unit cell size $\displaystyle\varepsilon_2=\frac{1}{4}$.
\begin{figure}[!htb]
\centering
\begin{minipage}[c]{0.4\textwidth}
  \centering
  \includegraphics[width=50mm]{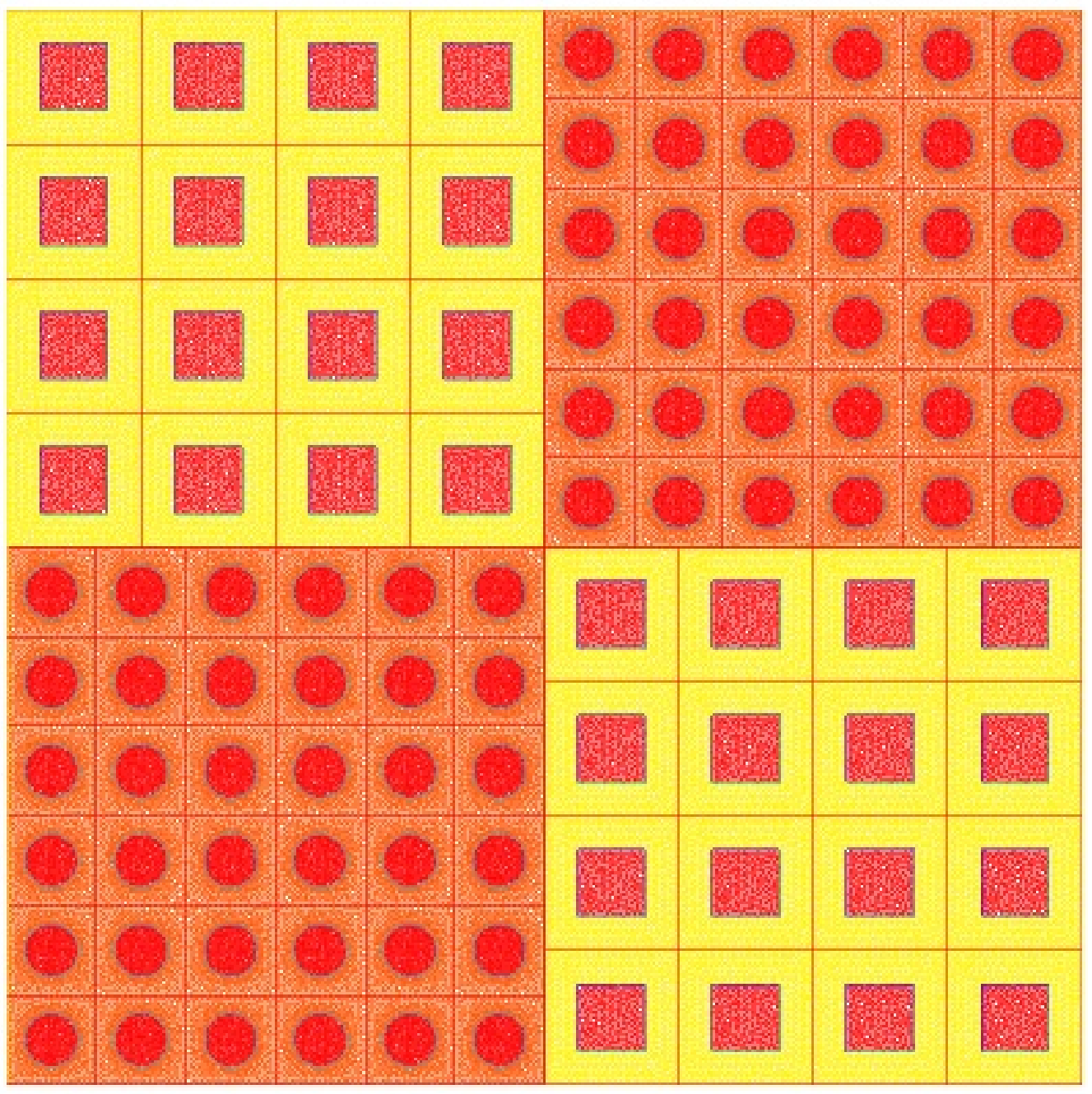}\\
  (a)
\end{minipage}
\begin{minipage}[c]{0.4\textwidth}
  \centering
  \includegraphics[width=50mm]{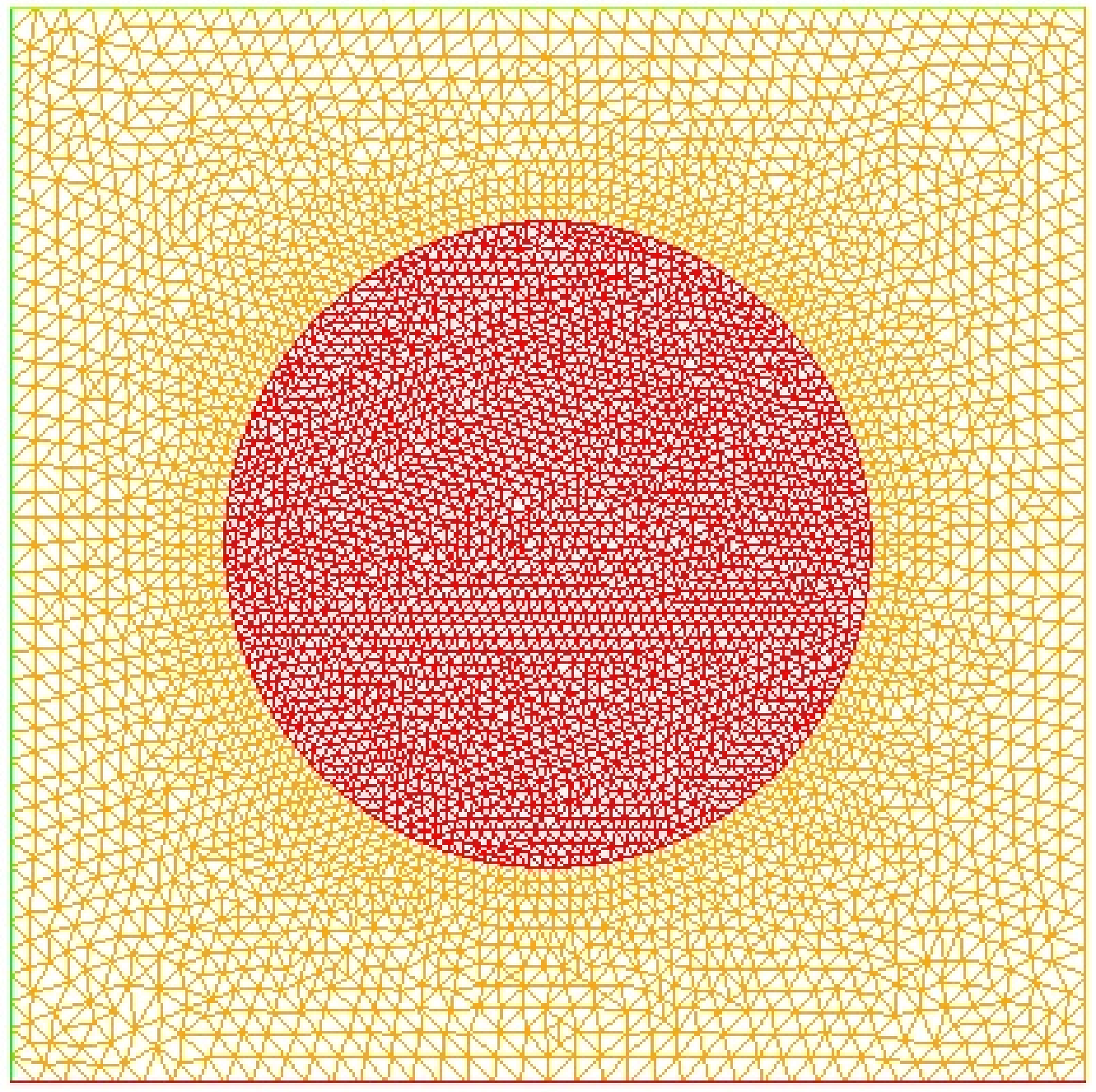}\\
  (b)
\end{minipage}
\begin{minipage}[c]{0.4\textwidth}
  \centering
  \includegraphics[width=50mm]{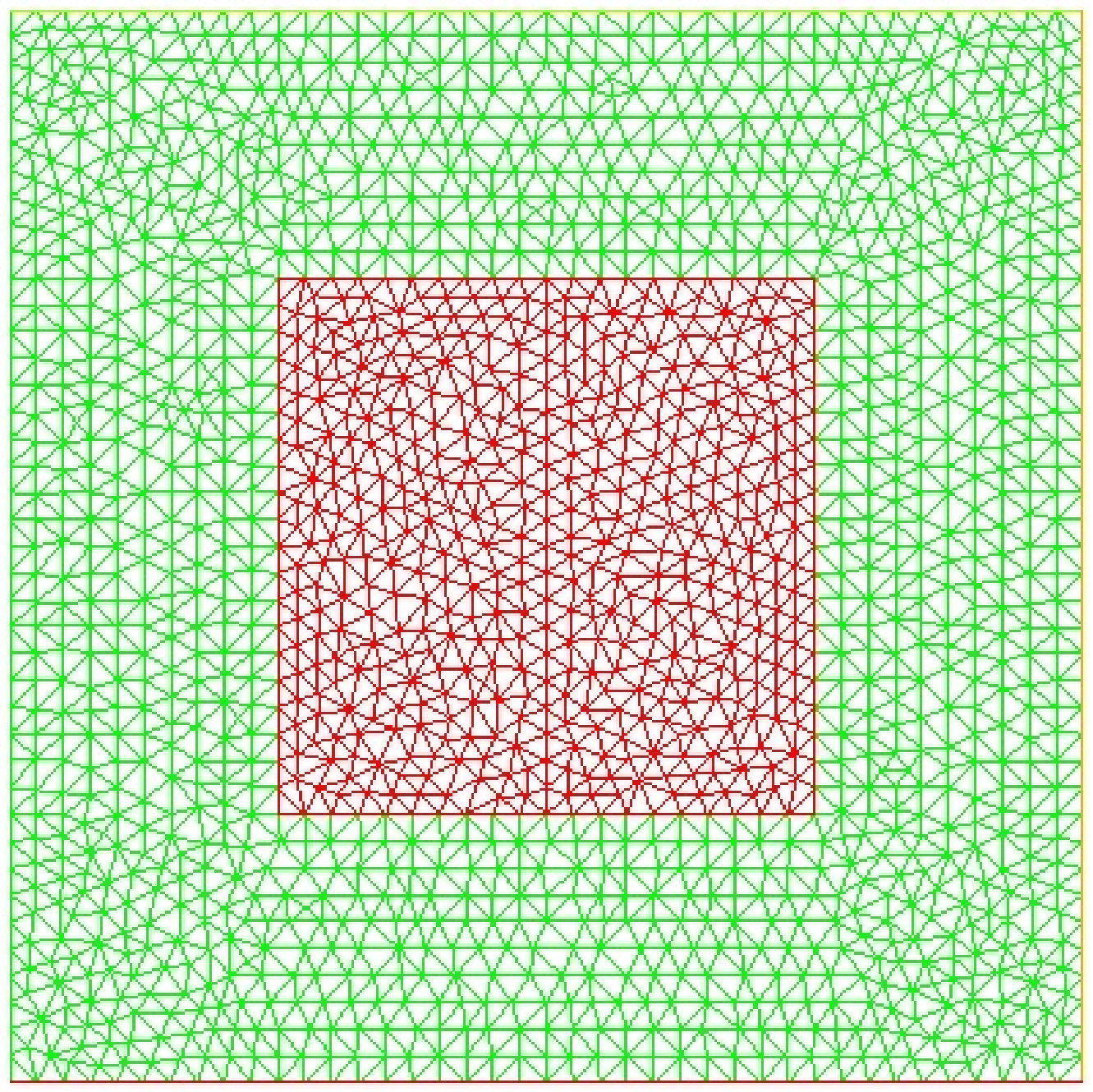}\\
  (c)
\end{minipage}
\caption{(a) The macrostructure domain $\Omega$; (b) the unit cell $Q^1$; (c) the unit cell $Q^2$.}\label{f8}
\end{figure}

In this example, the inclusion and matrix of $Q^1$ and $Q^2$ are defined as the identical materia. The detailed material property parameters are listed in Table 1.
\begin{table}[!htb]{\caption{Material property parameters}\label{t2}}
\centering
\begin{tabular}{ccc}
\hline
Property &Matrix of $\Omega_1\cup\Omega_3$ & Inclusion of $\Omega_1\cup\Omega_3$\\
\hline
Thermal conductivity$(W/(m\bm{\cdot} K))$ & 100.0 & 0.1 \\
\hline
Property &Matrix of $\Omega_2\cup\Omega_4$ & Inclusion of $\Omega_2\cup\Omega_4$\\
\hline
Thermal conductivity$(W/(m\bm{\cdot} K))$ & 100.0 & 0.1 \\
\hline
\end{tabular}
\end{table}

The data in problem (1) are given as follows:
\begin{equation}
h(\mathbf{x})=100J/(cm^3\bm{\cdot}s),\;\widehat{T}(\mathbf{x})=373.15 K\;\;\text{in}\;\;\partial\Omega
\end{equation}

Now, we implement the triangular mesh generation to multiscale problem (1), auxiliary cell problems and associated homogenized problem (12). The computational cost of FEM elements and nodes is listed in Table 2.
\begin{table}[!htb]{\caption{Comparison of computational cost}\label{t2}}
\centering
\begin{tabular}{ccc}
\hline
 & Original equation & Cell problem of $Q^1$\\
\hline
number of elements & 115216 & 3446\\
number of nodes    & 58473 & 1804\\
\hline
& Cell problem of $Q^2$ & Homogenized equation \\
\hline
number of elements & 3438 & 53286\\
number of nodes    & 1800 & 26944\\
\hline
\end{tabular}
\end{table}

Fig. 2 shows the numerical results for solutions $T^{\bm{\varepsilon}}$, $T^{(0)}$, $T^{(1\bm{\varepsilon})}$ and $T^{(2\bm{\varepsilon})}$, respectively.
\begin{figure}[!htb]
\centering
\begin{minipage}[c]{0.4\textwidth}
  \centering
  \includegraphics[width=50mm]{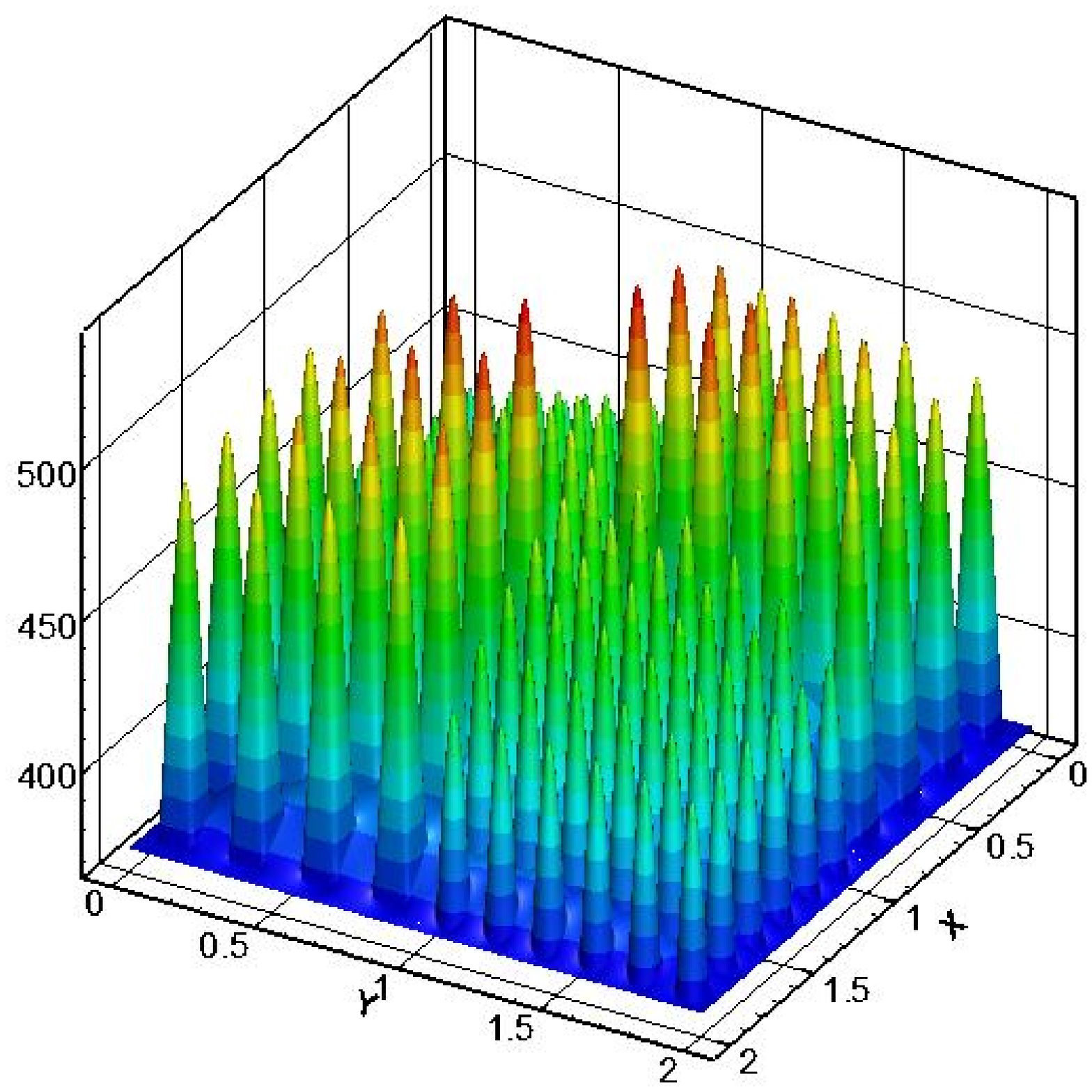}\\
  (a)
\end{minipage}
\begin{minipage}[c]{0.4\textwidth}
  \centering
  \includegraphics[width=50mm]{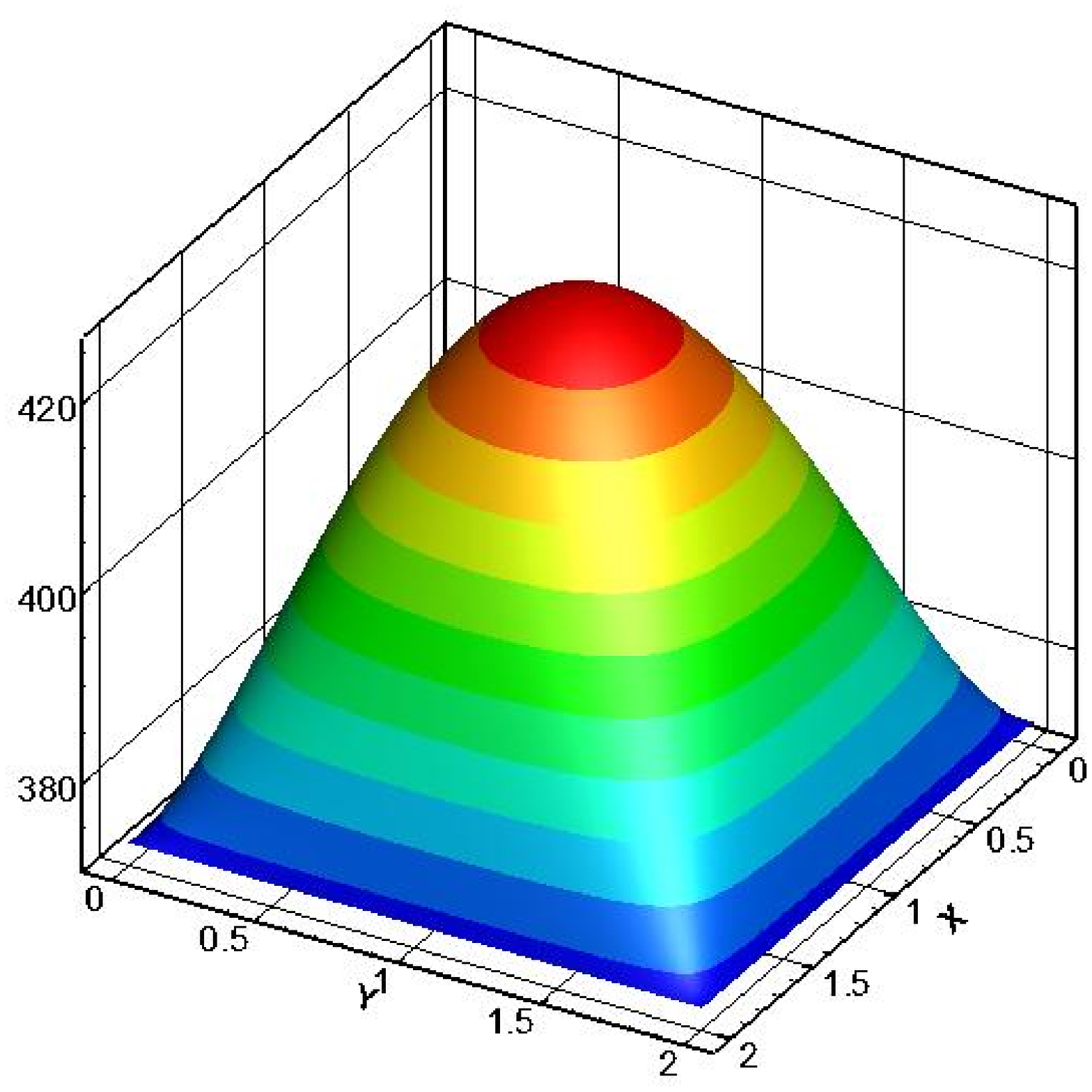}\\
  (b)
\end{minipage}
\begin{minipage}[c]{0.4\textwidth}
  \centering
  \includegraphics[width=50mm]{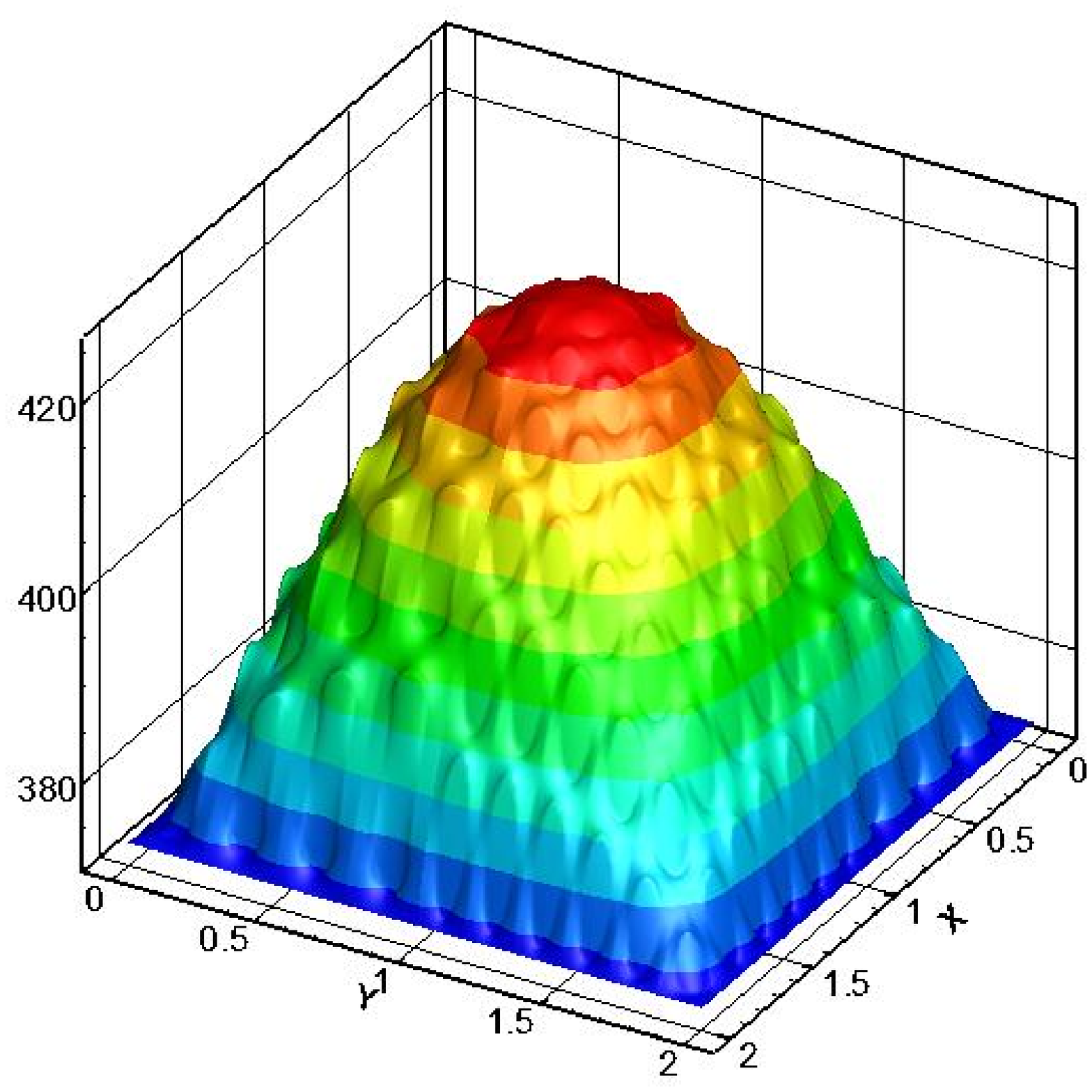}\\
  (c)
\end{minipage}
\begin{minipage}[c]{0.4\textwidth}
  \centering
  \includegraphics[width=50mm]{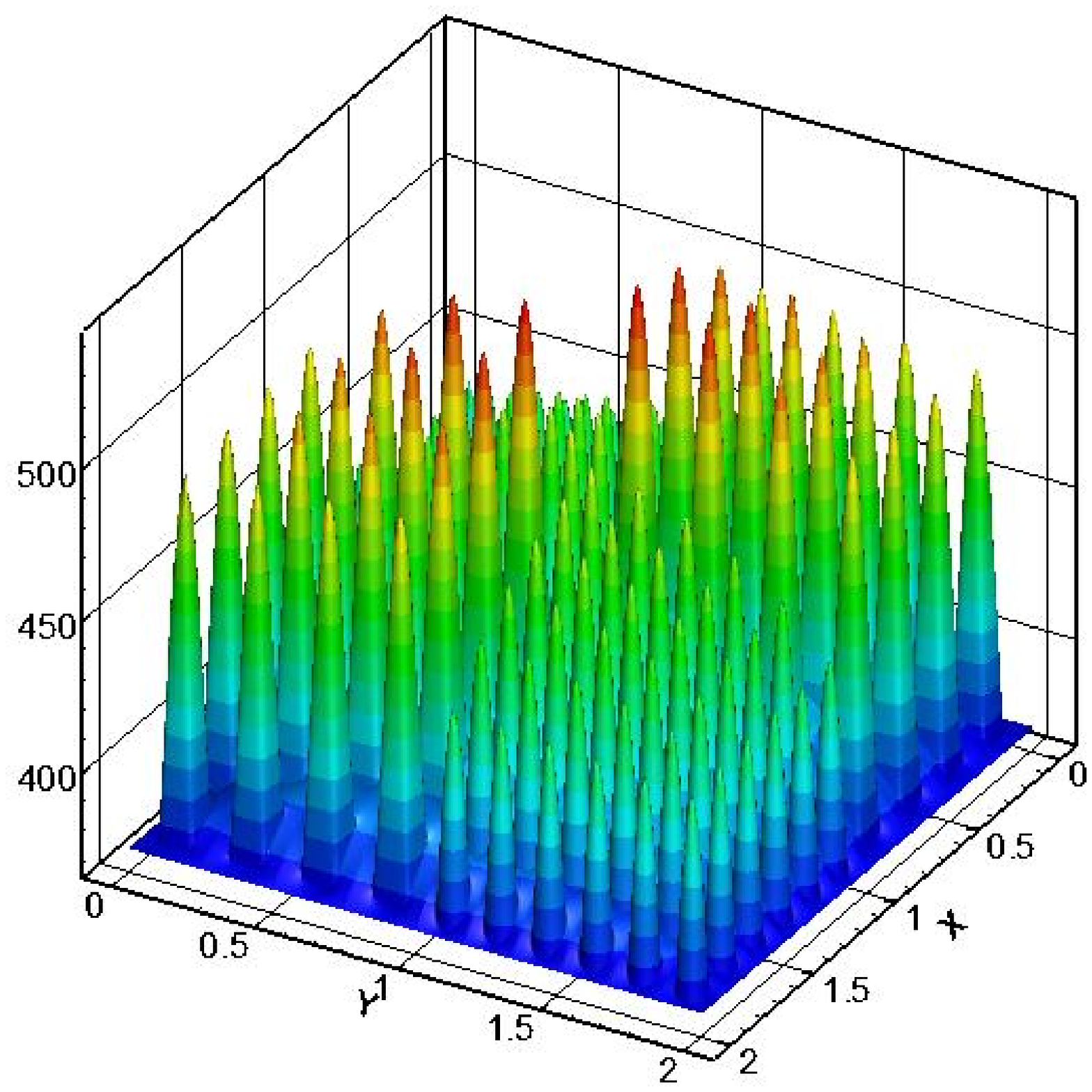}\\
  (d)
\end{minipage}
\caption{The temperature field: (a) $T^{\bm{\varepsilon}}$; (b) $T^{(0)}$; (c) $T^{(1\bm{\varepsilon})}$; (d) $T^{(2\bm{\varepsilon})}$.}\label{f8}
\end{figure}

After completing numerical computation, the relative $L^2$ norm error and $H^1$ semi-norm error of temperature field are listed in Table 3.
\begin{table}[h]{\caption{Comparison of relative errors}\label{t2}}
\centering
\begin{tabular}{cccc}
\hline
&Terror0 & Terror1 & Terror2\\
\hline
Percentage \% &6.3920 & 6.3951 & 0.0627\\
\hline
&TError0 & TError1 & TError2\\
\hline
Percentage \% &99.6400 & 99.3812 & 5.9648\\
\hline
\end{tabular}
\end{table}

From Table 2, one can see that the computational cost of SOTS method is much less than precise FEM. It means that the SOTS method can greatly save computer memory, which is very important in engineering computation. Fig. 2 demonstrates that only SOTS solution can accurately capture the micro-scale oscillating information due to heterogeneities in composites. From Table 3, we can conclude that only SOTS solution is almost the same as the precise FEM solution. By contrast, homogenized and FOTS solutions are far from enough to provide a high accuracy solution for multiscale problem (1).
\subsection{Example 2: different configurations with diverse material}
In this example, a 2D composite structures with two basic periodic configurations in four subdomains is considered. The macrostructure $\Omega$ is shown in Fig. 3, where $\Omega=(x,y)=[0,2]\times[0,2]cm^2$. Assume that $\Omega=\Omega_1\cup\Omega_2\cup\Omega_3\cup\Omega_4$ where $\Omega_1=[0,1]\times[0,1]cm^2$, $\Omega_2=[1,2]\times[0,1]cm^2$, $\Omega_3=[1,2]\times[1,2]cm^2$ and $\Omega_4=[0,1]\times[1,2]cm^2$. Moreover, $\Omega_1$ and $\Omega_3$ have the same unit cell $Q^1$ with periodic unit cell size $\displaystyle\varepsilon_1=\frac{1}{7}$. $\Omega_2$ and $\Omega_4$ have the same unit cell $Q^2$ with periodic unit cell size $\displaystyle\varepsilon_2=\frac{1}{5}$.
\begin{figure}[!htb]
\centering
\begin{minipage}[c]{0.4\textwidth}
  \centering
  \includegraphics[width=50mm]{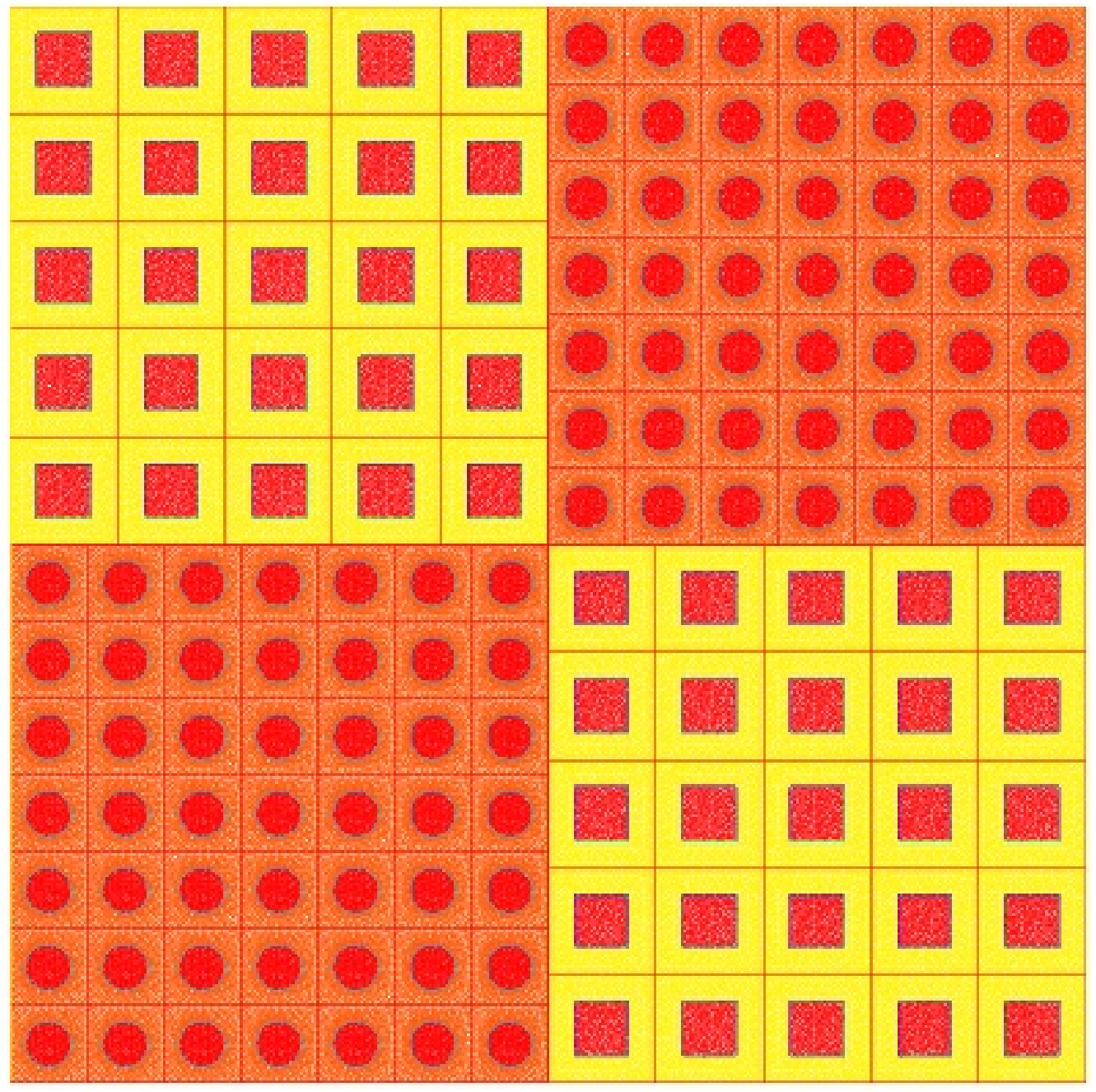}\\
  (a)
\end{minipage}
\begin{minipage}[c]{0.4\textwidth}
  \centering
  \includegraphics[width=50mm]{10.eps}\\
  (b)
\end{minipage}
\begin{minipage}[c]{0.4\textwidth}
  \centering
  \includegraphics[width=50mm]{11.eps}\\
  (c)
\end{minipage}
\caption{(a) The macrostructure domain $\Omega$; (b) the unit cell $Q^1$; (c) the unit cell $Q^2$.}\label{f8}
\end{figure}

In this example, the inclusion of $Q^1$ and $Q^2$ are defined as different materials. The detailed material property parameters are listed in Table 4.
\begin{table}[!htb]{\caption{Material property parameters}\label{t2}}
\centering
\begin{tabular}{ccc}
\hline
Property &Matrix of $\Omega_1\cup\Omega_3$ & Inclusion of $\Omega_1\cup\Omega_3$\\
\hline
Thermal conductivity$(W/(m\bm{\cdot} K))$ & 100.0 & 0.5 \\
\hline
Property &Matrix of $\Omega_2\cup\Omega_4$ & Inclusion of $\Omega_2\cup\Omega_4$\\
\hline
Thermal conductivity$(W/(m\bm{\cdot} K))$ & 100.0 & 0.1 \\
\hline
\end{tabular}
\end{table}

The data in problem (1) are given as follows:
\begin{equation}
h(\mathbf{x})=200J/(cm^3\bm{\cdot}s),\;\widehat{T}(\mathbf{x})=373.15 K\;\;\text{in}\;\;\partial\Omega
\end{equation}

Now, we implement the triangular mesh generation to multiscale problem (1), auxiliary cell problems and associated homogenized problem (12). The computational cost of FEM elements and nodes is listed in Table 5.
\begin{table}[!htb]{\caption{Comparison of computational cost }\label{t2}}
\centering
\begin{tabular}{ccc}
\hline
 & Original equation & Cell problem of $Q^1$\\
\hline
number of elements & 169904 & 3446\\
number of nodes    & 85985 & 1804\\
\hline
& Cell problem of $Q^2$ & Homogenized equation \\
\hline
number of elements & 3438 & 53286\\
number of nodes    & 1800 & 26944\\
\hline
\end{tabular}
\end{table}

Fig. 4 shows the numerical results for solutions $T^{\bm{\varepsilon}}$, $T^{(0)}$, $T^{(1\bm{\varepsilon})}$ and $T^{(2\bm{\varepsilon})}$, respectively.
\begin{figure}[!htb]
\centering
\begin{minipage}[c]{0.4\textwidth}
  \centering
  \includegraphics[width=50mm]{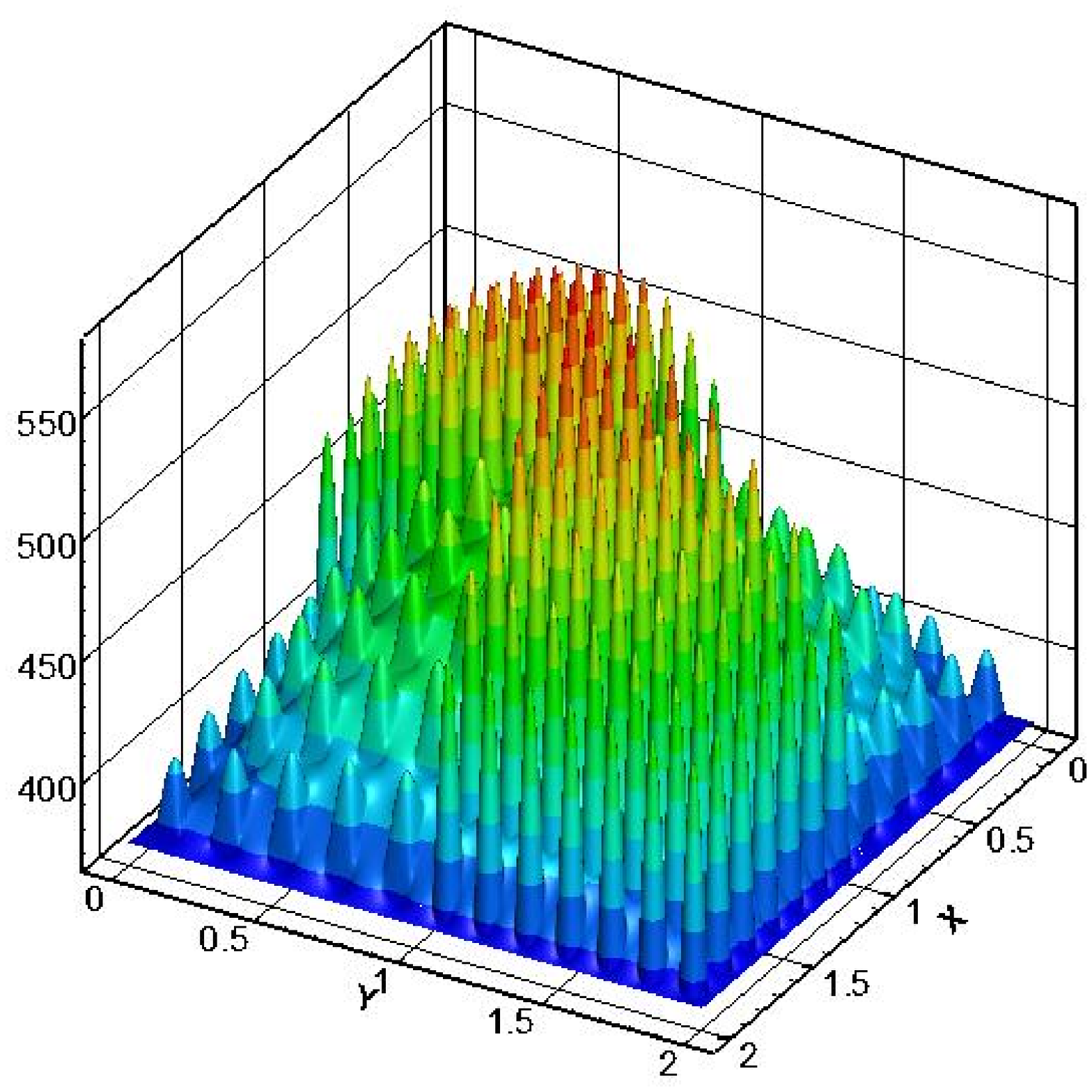}\\
  (a)
\end{minipage}
\begin{minipage}[c]{0.4\textwidth}
  \centering
  \includegraphics[width=50mm]{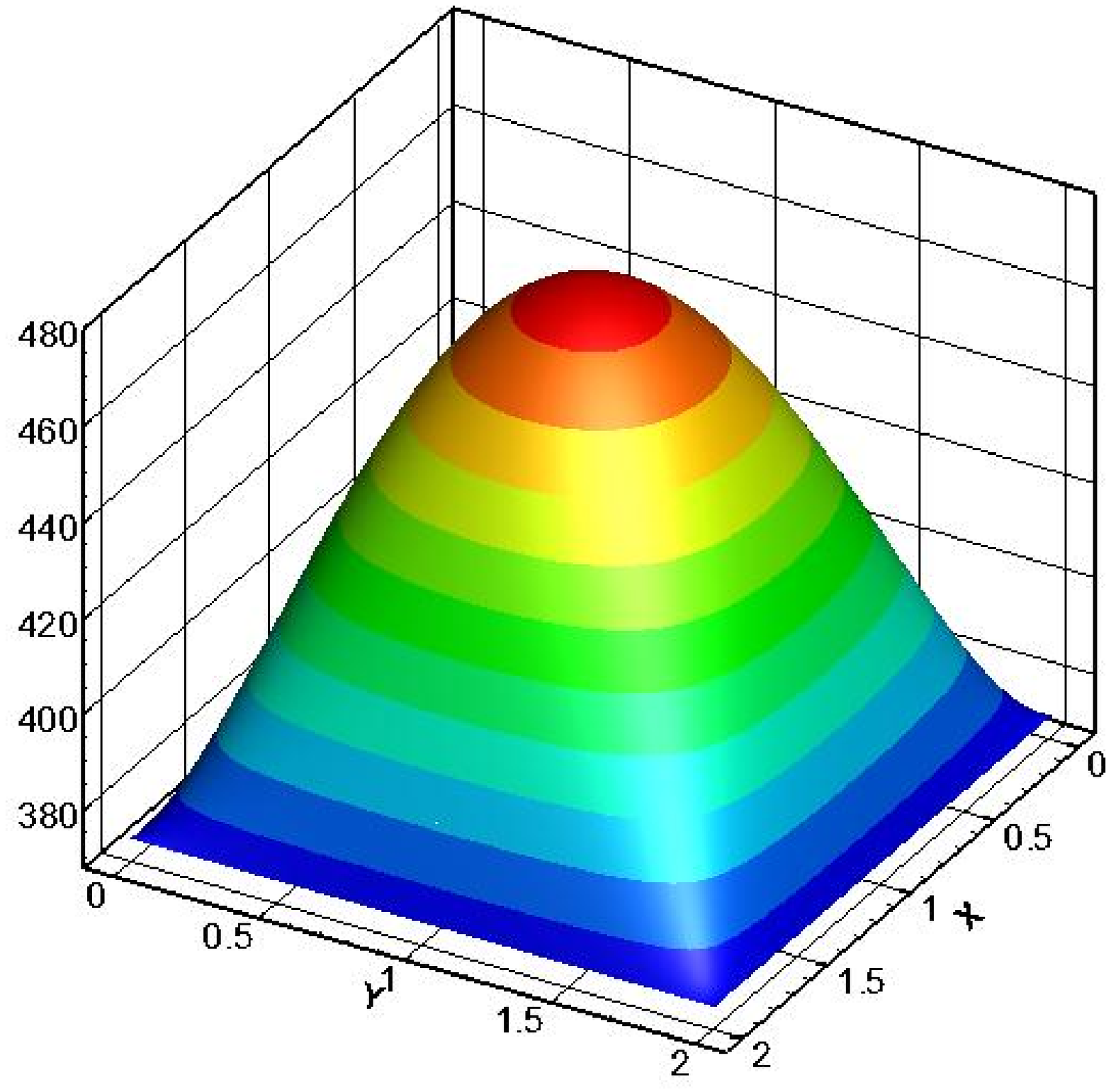}\\
  (b)
\end{minipage}
\begin{minipage}[c]{0.4\textwidth}
  \centering
  \includegraphics[width=50mm]{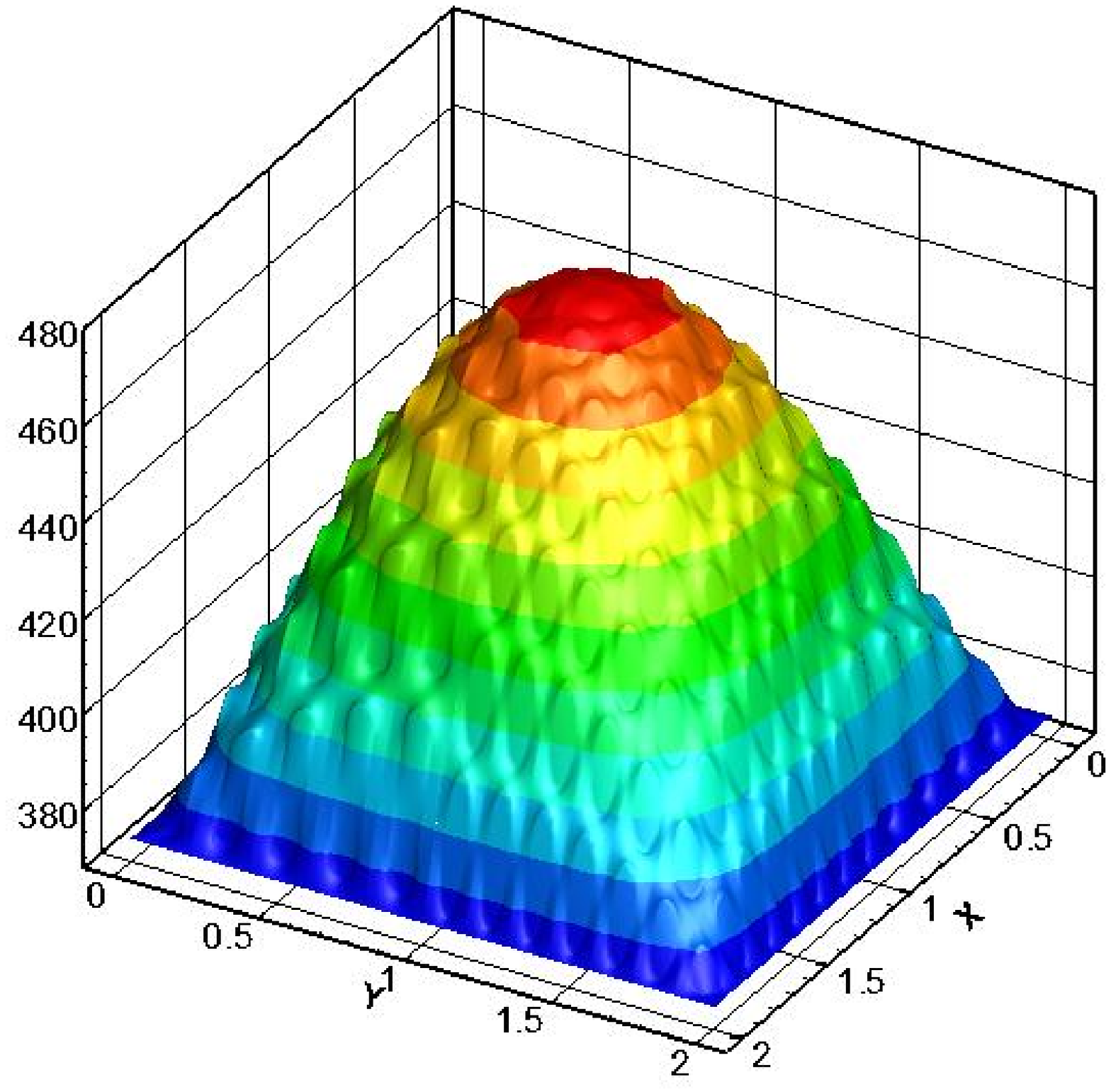}\\
  (c)
\end{minipage}
\begin{minipage}[c]{0.4\textwidth}
  \centering
  \includegraphics[width=50mm]{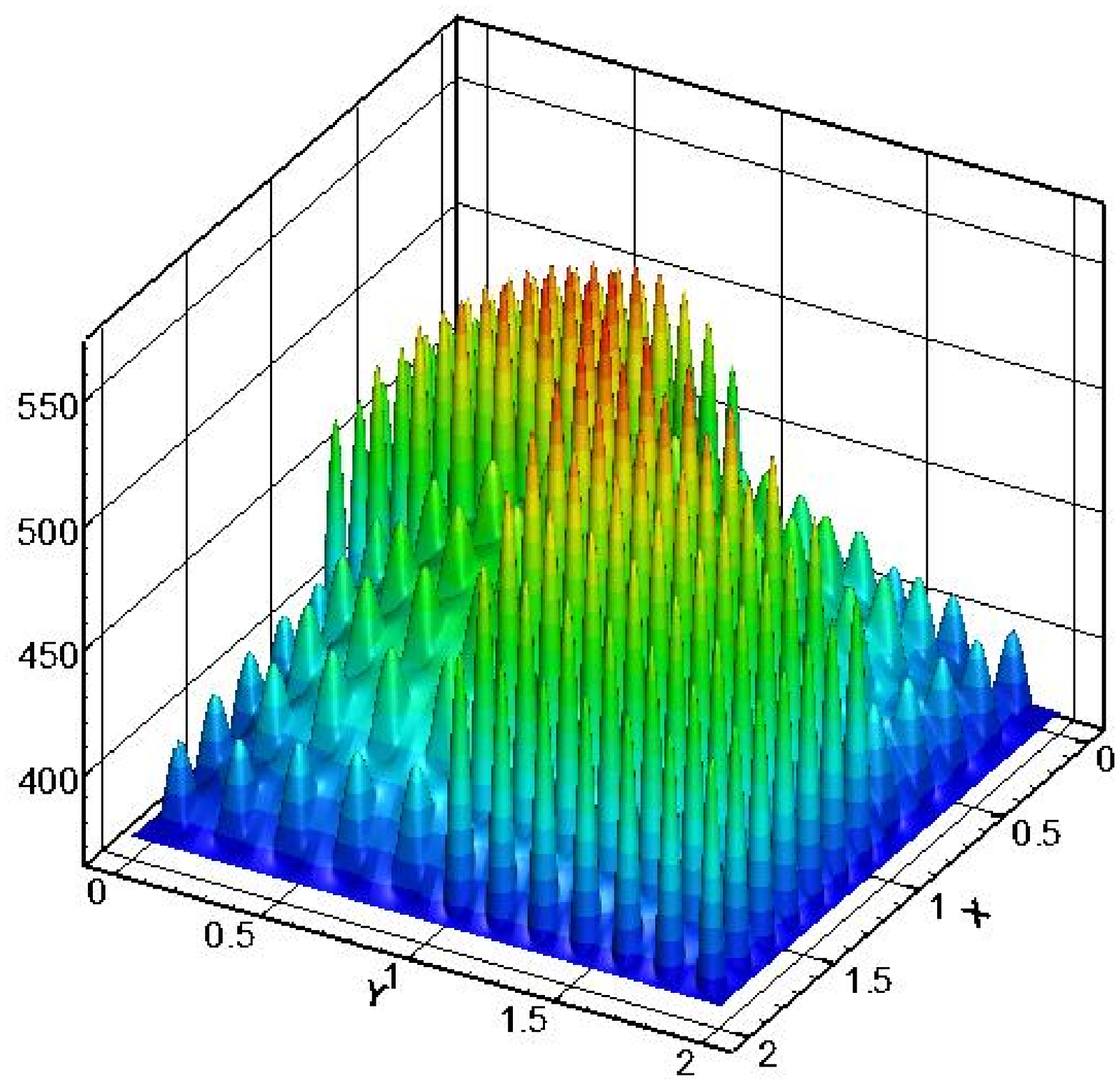}\\
  (d)
\end{minipage}
\caption{The temperature field: (a) $T^{\bm{\varepsilon}}$; (b) $T^{(0)}$; (c) $T^{(1\bm{\varepsilon})}$; (d) $T^{(2\bm{\varepsilon})}$.}\label{f8}
\end{figure}

Afterwards, the relative $L^2$ norm error and $H^1$ semi-norm error of temperature field are listed in Table 6.
\begin{table}[h]{\caption{Comparison of relative errors}\label{t2}}
\centering
\begin{tabular}{cccc}
\hline
&Terror0 & Terror1 & Terror2\\
\hline
Percentage \% &5.1486 & 5.1398 & 1.0740\\
\hline
&TError0 & TError1 & TError2\\
\hline
Percentage \% &99.0826 & 98.6436 & 8.8350\\
\hline
\end{tabular}
\end{table}

From Table 5, one can easily see that the computational cost of SOTS method still is much less than precise FEM. From Fig. 4 , it shows that the SOTS solution is much better than the homogenized and FOTS solutions for temperature field of multiscale problem (1). It is easy to see that only the SOTS solutions can provide enough numerical accuracy for engineering applications from Table 6. The accuracy of homogenized solutions and FOTS solutions is far from enough especially in the $H^1$ semi-norm sense.
\section{Conclusions}
In this paper, we develop a SOTS computational method for heat conduction problems of composite structures with diverse periodic configurations in different subdomains. The new contributions of this paper are the SOTS analysis, the error analysis in the pointwise sense and integral sense for the SOTS solutions, and associated SOTS numerical algorithm. Numerical experiments show that the SOTS numerical method we proposed is effective for multiscale problem (1). Furthermore, numerical results show that only SOTS solutions can accurately capture the microscale oscillating information and provide enough numerical accuracy for engineering applications, which support the theoretical results of this paper.
\section*{Acknowledgments}
This research was funded by the National Natural Science Foundation of China (No.11471262 and 11501449), the National Basic Research Program of China (No.2012CB025904), the State Scholarship Fund of China Scholarship Council (File No. 201606290191), and also supported by the Center for high performance computing of Northwestern Polytechnical University.




\bibliographystyle{model1-num-names}
\bibliography{paper}

\begin{thebibliography}{16}
\expandafter\ifx\csname natexlab\endcsname\relax\def\natexlab#1{#1}\fi
\providecommand{\url}[1]{\texttt{#1}}
\providecommand{\href}[2]{#2}
\providecommand{\path}[1]{#1}
\providecommand{\DOIprefix}{doi:}
\providecommand{\ArXivprefix}{arXiv:}
\providecommand{\URLprefix}{URL: }
\providecommand{\Pubmedprefix}{pmid:}
\providecommand{\doi}[1]{\href{http://dx.doi.org/#1}{\path{#1}}}
\providecommand{\Pubmed}[1]{\href{pmid:#1}{\path{#1}}}
\providecommand{\bibinfo}[2]{#2}
\ifx\xfnm\relax \def\xfnm[#1]{\unskip,\space#1}\fi
\bibitem[{Matine et~al.(2015)Matine, Boyard, Legrain, Jarny, and Cartraud}]{R1}
\bibinfo{author}{A.~Matine}, \bibinfo{author}{N.~Boyard},
  \bibinfo{author}{G.~Legrain}, \bibinfo{author}{Y.~Jarny},
  \bibinfo{author}{P.~Cartraud},
\newblock \bibinfo{title}{{Transient heat conduction within periodic
  heterogeneous media: A space-time homogenization approach}},
\newblock \bibinfo{journal}{International Journal of Thermal Sciences}
  \bibinfo{volume}{92} (\bibinfo{year}{2015}) \bibinfo{pages}{217--229}.
\bibitem[{Ma et~al.(2016)Ma, Cui, Li, and Wang}]{R2}
\bibinfo{author}{Q.~Ma}, \bibinfo{author}{J.~Cui}, \bibinfo{author}{Z.~Li},
  \bibinfo{author}{Z.~Wang},
\newblock \bibinfo{title}{{Second-order asymptotic algorithm for heat
  conduction problems of periodic composite materials in curvilinear
  coordinates}},
\newblock \bibinfo{journal}{Journal of Computational and Applied Mathematics}
  \bibinfo{volume}{306} (\bibinfo{year}{2016}) \bibinfo{pages}{87--115}.
\bibitem[{Yanga et~al.(2015)Yanga, Cui, Sun, and Ge}]{R3}
\bibinfo{author}{Z.~Yanga}, \bibinfo{author}{J.~Cui}, \bibinfo{author}{Y.~Sun},
  \bibinfo{author}{J.~Ge},
\newblock \bibinfo{title}{{Multiscale computation for transient heat conduction
  problem with radiation boundary condition in porous materials}},
\newblock \bibinfo{journal}{Finite Elements in Analysis and Design}
  \bibinfo{volume}{102-103} (\bibinfo{year}{2015}) \bibinfo{pages}{7--18}.
\bibitem[{Yang et~al.(2016)Yang, Cui, Wang, and Zhang}]{R4}
\bibinfo{author}{Z.~Yang}, \bibinfo{author}{J.~Cui}, \bibinfo{author}{Z.~Wang},
  \bibinfo{author}{Y.~Zhang},
\newblock \bibinfo{title}{{Multiscale computational method for nonstationary
  integrated heat transfer problem in periodic porous materials}},
\newblock \bibinfo{journal}{Numerical Methods for Partial Differential
  Equations} \bibinfo{volume}{32} (\bibinfo{year}{2016})
  \bibinfo{pages}{510--530}.
\bibitem[{Su et~al.(2011)Su, Xu, Dong, and Li}]{R5}
\bibinfo{author}{F.~Su}, \bibinfo{author}{Z.~Xu}, \bibinfo{author}{Q.~Dong},
  \bibinfo{author}{H.~Li},
\newblock \bibinfo{title}{{A second-order and two-scale computation method for
  heat conduction equation with rapidly oscillatory coefficients}},
\newblock \bibinfo{journal}{Finite Elements in Analysis and Design}
  \bibinfo{volume}{47} (\bibinfo{year}{2011}) \bibinfo{pages}{276--280}.
\bibitem[{Allegretto et~al.(2008)Allegretto, Cao, and Lin}]{R6}
\bibinfo{author}{W.~Allegretto}, \bibinfo{author}{L.~Cao},
  \bibinfo{author}{Y.~Lin},
\newblock \bibinfo{title}{{Multiscale asymptotic expansion for second order
  parabolic equations with rapidly oscillating coefficients}},
\newblock \bibinfo{journal}{Discrete and Continuous Dynamical Systems A}
  \bibinfo{volume}{20} (\bibinfo{year}{2008}) \bibinfo{pages}{543--576}.
\bibitem[{Wu et~al.(2014)Wu, Nie, and Yang}]{R7}
\bibinfo{author}{Y.~T. Wu}, \bibinfo{author}{Y.~F. Nie}, \bibinfo{author}{Z.~H.
  Yang},
\newblock \bibinfo{title}{{Comparison of four multiscale methods for elliptic
  problems}},
\newblock \bibinfo{journal}{CMES-Computer Modeling in Engineering $\&$
  Sciences} \bibinfo{volume}{99} (\bibinfo{year}{2014})
  \bibinfo{pages}{297--325}.
\bibitem[{Dong et~al.(2016)Dong, Nie, Yang, and Wu}]{R8}
\bibinfo{author}{H.~Dong}, \bibinfo{author}{Y.~Nie}, \bibinfo{author}{Z.~Yang},
  \bibinfo{author}{Y.~Wu},
\newblock \bibinfo{title}{{The Numerical Accuracy Analysis of Asymptotic
  Homogenization Method and Multiscale Finite Element Method for Periodic
  Composite Materials}},
\newblock \bibinfo{journal}{CMES-Computer Modeling in Engineering $\&$
  Sciences} \bibinfo{volume}{111} (\bibinfo{year}{2016})
  \bibinfo{pages}{395--419}.
\bibitem[{Cui(2001)}]{R9}
\bibinfo{author}{J.~Cui},
\newblock \bibinfo{title}{{Multiscale computational method for unified design
  of structure, components and their materials}},
\newblock in: \bibinfo{booktitle}{Proceedings on Computational Mechanics in
  Science and Engineering, CCCM-2001, Guangzhou, 5-8 December},
  \bibinfo{publisher}{Peking University Press}, \bibinfo{year}{2001}, pp.
  \bibinfo{pages}{33--43}.
\bibitem[{Cui(1996)}]{R10}
\bibinfo{author}{J.~Cui},
\newblock \bibinfo{title}{{The two-scale expression of the solution for the
  structure with several sub-domains of small periodic configurations}},
\newblock in: \bibinfo{booktitle}{WORKSHOP ON SCIENTIFIC COMPUTING 99, Hong
  Kong, 27-30 June}, \bibinfo{year}{1996}.
\bibitem[{Oleinik et~al.(1992)Oleinik, Shamaev, and Yosifian}]{R11}
\bibinfo{author}{O.~A. Oleinik}, \bibinfo{author}{A.~S. Shamaev},
  \bibinfo{author}{G.~A. Yosifian}, \bibinfo{title}{{Mathematical Problems in
  Elasticity and Homogenization}}, \bibinfo{publisher}{North-Holland:
  Amsterdam}, \bibinfo{year}{{1992}}.
\bibitem[{Cioranescu and Donato(1999)}]{R12}
\bibinfo{author}{D.~Cioranescu}, \bibinfo{author}{P.~Donato},
  \bibinfo{title}{{An Introduction to Homogenization}},
  \bibinfo{publisher}{Oxford University Press}, \bibinfo{year}{{1999}}.
\bibitem[{Wang et~al.(2015)Wang, Cao, and Wong}]{R13}
\bibinfo{author}{X.~Wang}, \bibinfo{author}{L.~Cao}, \bibinfo{author}{Y.~Wong},
\newblock \bibinfo{title}{{Multiscale computation and convergence for coupled
  thermoelastic system in composite materials}},
\newblock \bibinfo{journal}{{SIAM Multiscale Modeling and Simulation}}
  \bibinfo{volume}{{13}} (\bibinfo{year}{{2015}}) \bibinfo{pages}{{661--690}}.
\bibitem[{Cao(2006)}]{R14}
\bibinfo{author}{L.~Cao},
\newblock \bibinfo{title}{{Multiscale asymptotic expansion and finite element
  methods for the mixed boundary value problems of second order elliptic
  equation in perforated domains}},
\newblock \bibinfo{journal}{{Numerische Mathematik}} \bibinfo{volume}{{103}}
  (\bibinfo{year}{{2006}}) \bibinfo{pages}{{11--45}}.
\bibitem[{Dong and Cao(2009)}]{R15}
\bibinfo{author}{Q.~Dong}, \bibinfo{author}{L.~Cao},
\newblock \bibinfo{title}{{Multiscale asymptotic expansions methods and
  numerical algorithms for the wave equations of second order with rapidly
  oscillating coefficients}},
\newblock \bibinfo{journal}{{Applied Numerical Mathematics}}
  \bibinfo{volume}{{59}} (\bibinfo{year}{{2009}}) \bibinfo{pages}{3008--3032}.
\bibitem[{Lin and Zhu(1994)}]{R16}
\bibinfo{author}{Q.~Lin}, \bibinfo{author}{Q.~Zhu}, \bibinfo{title}{{The
  Preprocessing snd Preprocessing for the Finite Element Method}},
  \bibinfo{publisher}{Shanghai Scientific $\&$ Technical Publishers},
  \bibinfo{year}{{1994}}.

\end{thebibliography}







\end{document}